\input ./style/arxiv-general.cfg
\documentclass[aop,MSNbibl,seceqn,dvips]{arximspdf}
\makeatletter
   \@ifpackageloaded{graphicx}{}{\usepackage{graphicx}}
\makeatother

\doi{10.1214/15-AOP1002}
\volume{44}
\issue{2}
\pubyear{2016}
\firstpage{1308}
\lastpage{1340}
\docsubty{FLA}

\makeatletter
\newcommand{\rrvert}{\vert}
\newcommand{\rrVert}{\Vert}
\newcommand{\llvert}{\vert}
\newcommand{\llVert}{\Vert}
\newproclaim{Exa}{Example}
\newcommand{\eqref}[1]{(\ref{#1})}
\newtheorem{Thm}{Theorem}[section]
\newtheorem{Prop}[Thm]{Proposition}
\newtheorem{Lem}[Thm]{Lemma}
\newtheorem{Cor}[Thm]{Corollary}
\newtheorem{Ex}[Thm]{Example}
\newproclaim{Rem}[Thm]{Remark}
\newtheorem{Fact}[Thm]{Fact}
\newproclaim{Def}{Definition}[section]
\newcommand{\g}[1]{\mathfrak{#1}}
\newproclaim{Rmq}[Thm]{Remark}
\newproclaim{Rmqs}[Thm]{Remarks}
\makeatother

\begin{document}
\begin{frontmatter}

\title{Central limit theorem for linear groups}
\runtitle{CLT for linear groups}

\begin{aug}
\author[A]{\fnms{Yves}~\snm{Benoist}\ead[label=e1]{yves.benoist@math.u-psud.fr}}%
 \and
\author[B]{\fnms{Jean-Fran\c{c}ois}~\snm{Quint}\corref{}\ead[label=e2]{Jean-Francois.Quint@math.u-bordeaux1.fr}}
\runauthor{Y. Benoist and J.-F. Quint}
\affiliation{CNRS---Universit\'e Paris-Sud and CNRS---Universit\'e Bordeaux I}
\address[A]{D\'{e}partement de Math\'{e}matiques\\
CNRS---Universit\'e Paris-Sud\\
Batiment  425\\
91405 Orsay\\
 France\\
\printead{e1}}
\address[B]{D\'{e}partement de Math\'{e}matiques\\
CNRS---Universit\'e Bordeaux I\\
33405 Talence\\
France\\
\printead{e2}}
\end{aug}

%
\received{\smonth{4} \syear{2013}}
%
\revised{\smonth{6} \syear{2014}}

\begin{abstract}
We prove a central limit theorem
for random walks with finite variance on linear groups.
\end{abstract}

\begin{keyword}[class=AMS]
\kwd{22E40}
\kwd{60G42}
\kwd{60G50}
\end{keyword}
\begin{keyword}
\kwd{Random walk}
\kwd{martingale}
\kwd{stationary measure}
\kwd{cocycle}
\kwd{semisimple group}
\end{keyword}
\end{frontmatter}

\section{Introduction}
\label{secintro}

\subsection{Central limit theorem for linear groups}
\label{secintcltlin}

Let $V=\mathbb{R}^d$, $G=\operatorname{GL}(V)$ and $\mu$ be a Borel probability
measure on $G$. We fix a norm $\|\cdot\|$ on $V$.
For $n\geq1$, we denote by $\mu^{*n}$ the $n${th}-convolution power
$\mu*\cdots*\mu$.
We assume that the first moment $\int_G \log N( g)\,\mathrm{d}\mu
(g)$ is finite,
where $N(g)=\max(\| g\|,\| g^{-1}\|)$.
We denote by $\lambda_1$ the first Lyapunov exponent of $\mu$, that is,
%
\begin{equation}
\label{eqnlalogmun} \lambda_1 :=\lim_{n\rightarrow\infty}
\frac{1}{n}\int_G\log\| g\|\,\mathrm{d}
\mu^{*n}(g).
\end{equation}
Let $g_1,\ldots, g_n,\ldots$ be random elements of $G$
chosen independently with law $\mu$.
The Furstenberg law of large numbers
describes the behavior of the random variables
$\log\| g_n\cdots g_1\|$.
It states that, almost surely,
%
\begin{equation}
\label{eqnlaloggn} \lim_{n\rightarrow\infty}\frac{1}n\log\|
g_n\cdots g_1\| =\lambda_1.
\end{equation}
In this paper, we will prove that, under suitable conditions, the
variables 
$\log\| g_n\cdots g_1\|$ satisfy a central limit theorem (CLT), that
is, the renormalized variables
$\frac{\log\| g_n\cdots g_1\|-n\lambda_1}{\sqrt{n}}$ converge in
law to a
nondegenerate Gaussian variable.

Let $\Gamma_\mu$ be the semigroup spanned by the support of $\mu$.
We say that
$\Gamma_\mu$ acts strongly irreducibly on $V$
if no proper finite union of vector subspaces of $V$ is $\Gamma_\mu$-invariant.

\begin{Thm}
\label{thcltmat}
Let $V=\mathbb{R}^d$, $G=\operatorname{GL}(V)$ and
$\mu$ be a Borel probability measure on $G$
such that $\Gamma_\mu$ has unbounded image in $\operatorname{PGL}(V)$,
$\Gamma_\mu$ acts strongly irreducibly on $V$,
and the second moment
$\int_G (\log N( g))^2\,\mathrm{d}\mu(g)$ is finite.
Let $\lambda_1$ be the first Lyapunov exponent of $\mu$.
Then there exists
$ \Phi>0$ such that,
for any bounded continuous function $F$ on $\mathbb{R}$, one has
%
\begin{equation}
\label{eqncltmat} \lim_{n\rightarrow\infty} \int_GF
\biggl(\frac{\log\| g\|-n\lambda_1}{\sqrt{n}} \biggr)\, \mathrm{d}\mu^{*n}(g) = \int
_{\mathbb{R}}F(s)\frac{e^{-{s^2}/({2\Phi})}}{\sqrt{2\pi
\Phi}}\,\mathrm{d}s.
\end{equation}
\end{Thm}

\begin{Rmqs}
We will see that under the same assumptions the variables
$\log\| g_n\cdots g_1\|$ also satisfy a law of the iterated logarithm
(LIL), that is, almost surely, the set of cluster points of the sequence
$\frac{\log\| g_n\cdots g_1\|-n\lambda_1}{\sqrt{2\Phi\, n\log\log
n}}$ is
equal to the interval $[-1,1]$.

According to a result of Furstenberg,
when moreover $\Gamma_\mu$ is included in the group $\operatorname{SL}(V)$,
the first Lyapunov exponent is positive: $\lambda_1>0$.

For every nonzero $v$ in $V$ and $f$ in $V^*$,
the variables $\log\| g_n\cdots g_1v\|$ and
$\log|f( g_n\cdots g_1v)|$ also satisfy the CLT and the LIL.

Such a central limit theorem is not always
true when the action of $\Gamma_\mu$ is only assumed to be irreducible:
in this case the variables $\frac{\log\| g_n\cdots g_1\|-n\lambda
_1}{\sqrt{n}}$
still converge in law but the limit is not always a Gaussian variable
(see Example~\ref{exanongau}).

We will deduce easily a multidimensional version of
this CLT (Theorem~\ref{thcltmatbis})
and interpret it as a CLT
for real semisimple groups (Theorem~\ref{thcltmatter}),
generalizing Goldsheid and Guivarc'h CLT in \cite{GoGu}.
Most of our results are true over any local field $\mathbb{K}$
with no changes in the proofs.
\end{Rmqs}

\subsection{Previous results}
\label{secpreres}
Let us give a historical perspective about this theorem.
The existence of such a ``noncommutative CLT''
was first guessed by Bellman in \cite{Bel54}.
Such a theorem has first been proved by Furstenberg and Kesten in \cite{FuKe}
for semigroups of positive matrices under an $\mathrm{L}^{2+\varepsilon}$
assumption
for some $\varepsilon>0$.
It was then extended by Le Page
in \cite{LeP} for more general semigroups when the law has a finite
exponential moment,
that is, when
there exists $\alpha>0$ such that $\int_GN(g)^\alpha\,\mathrm{d}\mu
(g)<\infty$.
Thanks to later works of Guivarc'h and Raugi in \cite{GR85} and
Gol'dshe\u{\i}d and Margulis in \cite{GoMa}, the assumptions in the Le Page theorem
were clarified: the sole remaining but still unwanted assumption was
that $\mu$ had a finite exponential moment.

Hence, the purpose of our Theorem~\ref{thcltmat} is to replace this finite exponential moment assumption
by a finite second moment assumption.
Such a finite second moment assumption is optimal.

Partial results have been obtained recently in this direction.
Tutubalin in \cite{Tu77} has proved Theorem~\ref{thcltmat}
when the law $\mu$ is assumed to have a density.
Jan in his thesis (see \cite{Jan00}) has extended the Le Page theorem
under the assumption that all the $p$-moments of $\mu$ are finite.
Hennion in \cite{Hen97} has proved Theorem~\ref{thcltmat}
in the case of semigroups of positive matrices.

There exist a few books and surveys (\cite{BoLa}, \cite{Fur} or \cite
{BQRW}) about
this theory of ``products of random matrices.'' This theory has had recently
nice applications to the study of discrete subgroups of Lie groups
(as in \cite{G90,BFLM} or \cite{BQI}).
These applications motivated our interest
in a better understanding of this CLT.

\subsection{Other Central Limit Theorems}
\label{seccltoth}

The method we introduce in this paper is very flexible
since it does not rely on a spectral gap property.
In the forthcoming paper \cite{BQHG}, we will adapt this method to
prove the CLT in other
situations where the CLT is only known
under a finite exponential moment assumption:
\begin{longlist}[--]
\item[--] The CLT for free groups due to Sawyer--Steger in \cite{SS87}
and Ledrappier in \cite{Led01}.
\item[--] The CLT for Gromov hyperbolic groups due to Bjorklund in \cite{Bjo10}.
\end{longlist}

\subsection{Strategy}
\label{secintstr}
We explain now in few words the strategy of the proof
of our central limit Theorem~\ref{thcltmat}.
We want to prove the central limit theorem for the random variables
$\kappa(g_n\cdots g_1)$
where the quantity
%
\begin{equation}
\label{eqnkaglogng} \kappa(g):=\log\| g\|
\end{equation}
controls the size of the element $g$ in $G$.
Let $X:=\mathbb{P}(V)$ be the projective space of the vector space
$V:=\mathbb{R}^d$.
Since this function $\kappa$ on $G$ is closely related
to the ``norm cocycle'' $\sigma \dvtx  G\times X\rightarrow\mathbb{R}$
given by
%
\begin{equation}
\label{eqnsigxloggvv} \sigma(g,x):=\log\frac{\| gv\|}{\| v\|},
\end{equation}
for $g$ in $G$ and $x=\mathbb{R} v$ in $\mathbb{P}(V)$,
we are reduced to prove, for every $x$ in $X$,
a central limit theorem for the random variables
$\sigma(g_n\cdots g_1,x)$.

We will follow Gordin's method. This method has been introduced in
\cite
{Go69} and \cite{GoLi} and has been often used since then; see, for instance,
\cite{Leb99,Bjo10}. See also \cite{Der06}
and \cite{BoGo}, Appendix, for a survey of this method and
\cite{BoGo}, Section~2.4, for the use of this method in order to prove
a CLT and an invariance principle in the context
of products of independent random matrices.

Following Gordin's method means that, we will replace, adding a
suitable coboundary, this cocycle $\sigma$ by
another cocycle
$\sigma_0$ for which the ``expected increase'' is constant,
that is, such that
\[
\int_G\sigma_0(g,x)\,\mathrm{d}\mu(g)=
\lambda_1
\]
for all $x$ in $X$.
This will allow us to use the classical central limit theorem for martingales
due to Brown in \cite{Bro71}.
In order to find this cocycle $\sigma_0$, we have to
find a continuous function $\psi\in\mathcal{C}^0(X)$ which satisfies the
following cohomological equation
%
\begin{equation}
\label{eqnppsmps} \varphi=\psi-P_\mu\psi+\lambda_1,
\end{equation}
where $P_\mu\psi$ is the averaged function
\[
P_\mu\psi \dvtx x\mapsto\int_G\psi(gx)\,
\mathrm{d}\mu(g)
\]
and where $\varphi\in\mathcal{C}^0(X)$ is the expected increase of
the cocycle
$\sigma$
%
\begin{equation}
\label{eqnphxintsigx} \varphi \dvtx x\mapsto\int_G\sigma(g,x)\,
\mathrm{d}\mu(g).
\end{equation}

The classical strategy to solve this cohomological equation relies
on spectral properties of this operator $P_\mu$.
These spectral properties
might not be valid under a finite second moment assumption.
This is where our strategy differs from the classical strategy:
we solve this cohomological equation by giving an explicit formula
for the solution $\psi$ in terms of the $\check{\mu}$-stationary
measure $\nu^*$
on the dual projective space $\mathbb{P}(V^*)$, where $\check\mu$ is the
image of $\mu$
by $g\mapsto g^{-1}$. This formula is
%
\begin{equation}
\label{eqnpintlogde} \psi(x)=\int_{\mathbb{P}(V^*)}\log\delta(x,y)\,\mathrm{d}
\nu^*(y),
\end{equation}
where
$\delta(x,y)= \frac{|f(v)|}{\| f\|  \| v\|}$,
for $x=\mathbb{R} v$ in $\mathbb{P}(V)$ and $y=\mathbb{R}f$ in
$\mathbb{P}(V^*)$
(Proposition~\ref{prosolcohequ}).

The main issue is to check that this integral is finite,
that is, the stationary measure $\nu^*$ is $\log$-regular,
when the second moment of $\mu$ is finite (Proposition~\ref{prohausdimpv}).

Let us recall the Hsu--Robbins theorem which seems at a first glance unrelated.
This theorem is a strengthening of the classical law of large numbers
for centered square-integrable independent identically
distributed random real variables $(\varphi_n)_{n\geq1}$.
This theorem tells us that the averages $\frac{1}n(\varphi_1+\cdots
+\varphi_n)$
converge completely to $0$, that is,
for all $\varepsilon>0$, the following series converge:
%
\begin{equation}
\label{eqnsumpph1phn} \sum_{n\geq1}\mathbb{P}\biggl(
\frac{1}n|\varphi_1+\cdots+\varphi _n|>
\varepsilon\biggr)<\infty.
\end{equation}
The key point to prove the $\log$-regularity of the stationary measure
$\nu^*$
is to prove an analogue of the Hsu--Robbins theorem for martingales
under a suitable condition of domination by a square-integrable function
(Theorem~\ref{thhsurobmar})
and to deduce from it another analogue of the Hsu--Robbins theorem
for the Furstenberg law of large numbers (Proposition~\ref{procomconlin}).

Another important ingredient in the proof of the $\log$-regularity of
$\nu^*$
is the simplicity of the first Lyapunov
exponent due to Guivarc'h in \cite{Gui81} and \cite{GR85}.

\subsection{Plan}
\label{secintpla}

In Section~\ref{seclimmar}, we prove the complete convergence
in the law of large numbers for martingales with square-integrable increments
and we recall the central limit theorem for these martingales with
square-integrable increments.

In Section~\ref{seccencoc}, we prove a large deviations estimate in
the Breiman law of large numbers for functions over a Markov--Feller chain,
we deduce the complete convergence
in the law of large numbers for square-integrable cocycles over random walks
and the central limit theorem when the cocycle is centerable.

In Section~\ref{seccltlin}, we prove successively the complete convergence
in the Furstenberg law of large numbers, the log-regularity of the corresponding
stationary measure on the projective space, the centerability
of the norm cocycle and the central limit Theorem \ref{thcltmat}.
We end this chapter by the multidimensional version of this central
limit theorem.

\section{Limit theorems for martingales}
\label{seclimmar}
%
We collect in this chapter the limit theorems for martingales that we
will need
in Chapter~\ref{seccencoc}.

\subsection{Complete convergence for martingales}
\label{secllnmar}
%
In this section, we prove the complete convergence
in the law of large numbers for martingales.

Let $(\Omega,\mathcal{B},\mathbb{P})$ be a
probability space. We first recall that
a sequence $X_n$ of random variables
converges completely to $X_\infty$, if, for all $\varepsilon>0$,
$\sum_{n\geq1}\mathbb{P}(|X_n-X_\infty|\geq\varepsilon)<\infty$.
By the Borel--Cantelli lemma, complete convergence implies
almost sure convergence.
We recall now the following classical result due to Baum and Katz in
\cite{BaKa}.

\begin{Fact}
\label{thbaumkatz}
\textup{Let $p\geq1$, let $(\varphi_n)_{n\geq1}$ be
independent identically distributed real random variables and
$S_n=\varphi_1 +\cdots+\varphi_n$. The following statements are
equivalent:
\begin{longlist}[(ii)]
\item[(i)]  $\mathbb{E}|\varphi_1|^p<\infty$ and $\mathbb{E}(\varphi
_1)=0$,
\item[(ii)]  $\sum_{n\geq1}n^{p-2}\mathbb{P}(|S_n|\geq n\varepsilon
)<\infty $,
for all $\varepsilon>0$.
\end{longlist}}
\end{Fact}

When $p=2$ the implication $\textup{(i)}\Rightarrow\textup{(ii)}$ is due to Hsu--Robbins
\cite{HsRo} and
the converse is due to Erd\H os \cite{Er49}.
In this case, condition $\textup{(ii)}$ means that the sequence $\frac{1}{n} S_n$
converges completely toward $0$.

When $p=1$, this fact is due to Spitzer \cite{Spit}.

Our aim is to prove the following generalization of Baum--Katz theorem
to martingales.
Let $\mathcal{B}_0\subset\cdots\subset\mathcal{B}_n\subset\cdots$
be sub-$\sigma$-algebras of $\mathcal{B}$.
We recall that a \textit{martingale difference} is a sequence
$(\varphi_n)_{n\geq1}$ of integrable random variables
on $\Omega$ such that $\mathbb{E}(\varphi_n\vert \mathcal
{B}_{n-1})=0$ for all $n\geq1$.

\begin{Thm}
\label{thhsurobmar}
Let $p>1$, let $(\varphi_n)_{n\geq1}$ be a martingale difference
and $S_n:=\varphi_1+\cdots+\varphi_n$ the corresponding martingale.
We assume that there exists a positive function $\varphi$ in $\mathrm{L}^p(\Omega)$
such that, for $n\geq1$, $t>0$,
%
\begin{equation}
\label{eqnphnph} \mathbb{E}(\mathbf{1}_{\{|\varphi_n|>t\}}\vert \mathcal{B}_{n-1})
\leq \mathbb{P}\bigl(\{\varphi>t\}\bigr)\qquad\mbox{almost surely}.
\end{equation}
Then there exist constants $C_n=C_n(p,\varepsilon,\varphi)$ such that,
for $n\geq1$, $\varepsilon>0$,
%
\begin{equation}
\label{eqnpsncn} \mathbb{P}\bigl(|S_n|>n\varepsilon\bigr)\leq C_n
\quad\mbox{and}\quad \sum_{n\geq1}n^{p-2}C_n<
\infty.
\end{equation}
\end{Thm}

The fact that the constants $C_n$ are controlled by the dominating
function $\varphi$
will be important in our applications.
A related theorem was stated in \cite{Sto07} for $p>2$.
The extension to the case $p=2$ is crucial for our applications.
We stated our result for $p>1$ since the proof is not very different
when $p=2$.

\begin{pf*}{Proof of Theorem~\protect\ref{thhsurobmar}}
Our proof combines the original proof
of the Baum--Katz theorem with Burkholder inequality.
Since $p>1$, we pick $\gamma<1$ such that $\gamma>\frac{p+1}{2p}$.
We set, for $k\leq n$,
%
\begin{equation}
\label{eqnphnktn} \varphi_{n,k}:=\varphi_k
\mathbf{1}_{\{ |\varphi_k|\leq n^\gamma\}} \quad\mbox{and}\quad T_n:=\sum
_{1\leq k\leq n}\varphi_{n,k}.
\end{equation}
In order to lighten the calculations, we also set
%
\begin{equation}
\label{eqnolphnk} \overline\varphi_{n,k}:=\varphi_{n,k}-
\mathbb{E}(\varphi_{n,k}\vert \mathcal{B}_{k-1}) \quad
\mbox{and}\quad \overline{T}_n:=\sum_{1\leq k\leq n}
\overline\varphi_{n,k}
\end{equation}
so that, for all $n\geq1$,
the finite sequence $(\overline\varphi_{n,k})_{1\leq k\leq n}$ is
also a
difference martingale.
We can assume $\varepsilon=3$.
We will decompose the event
$A_n:=\{|S_n|>3n\}$ into four pieces
%
\begin{equation}
\label{eqnananan} A_n\subset A_{1,n}\cup A_{2,n}
\cup A_{3,n}\cup A_{4,n}.
\end{equation}
The events $A_{i,n}$ are given by
\begin{eqnarray*}
A_{1,n} &:=& \bigl\{ \mbox{there exists } k\leq n \mbox{
such that } |\varphi_k|>n\bigr\},
\\
A_{2,n} &:=& \bigl\{ \mbox{there exist } k_1<
k_2\leq n \mbox{ such that } |\varphi_{k_1}|>n^\gamma,
|\varphi_{k_2}|>n^\gamma\bigr\},
\\
A_{3,n} &:=& \bigl\{ |T_n-\overline T_n|> n\bigr\},
\\
A_{4,n} &:=& \bigl\{ |\overline T_n|> n\bigr\}.
\end{eqnarray*}
The inclusion \eqref{eqnananan} is satisfied since, when none of the
four events $A_{i,n}$ is satisfied, one has $|S_n|\leq3n$.
We will find, for each piece $A_{i,n}$, a constant
$C_{i,n}=C_{i,n}(p,\varepsilon,\varphi)$ such that
$\mathbb{P}(A_{i,n})\leq C_{i,n}$ and $\sum_{n\geq
1}n^{p-2}C_{i,n}<\infty$.

\begin{longlist}[{}]
\item[\textit{First piece}.] One computes, using the domination \eqref{eqnphnph},
\[
\mathbb{P}(A_{1,n})\leq C_{1,n}:=n \mathbb{P}(\varphi>n)
\]
and
\[
\sum_{n\geq1}n^{p-2}C_{1,n}= \sum
_{n\geq1}n^{p-1}\mathbb{P}(\varphi>n)\leq
\frac{1}p \mathbb {E}\bigl((\varphi+1)^p\bigr)
\]
which is finite since the dominating function $\varphi$ is $\mathrm{L}^p$-integrable.

\item[\textit{Second piece}.]
One computes, using the domination \eqref{eqnphnph},
\[
\mathbb{P}(A_{2,n})\leq C_{2,n}:=n^2 \mathbb{P}
\bigl(\varphi>n^\gamma\bigr)^2
\]
and, using Chebyshev's inequality,
\[
\sum_{n\geq1}n^{p-2}C_{2,n}\leq
\sum_{n\geq1}n^{p-2\gamma p}\bigl(\mathbb{E}\bigl(
\varphi^p\bigr)\bigr)^2
\]
which is finite since $\gamma>\frac{p+1}{2p}$.

\item[\textit{Third piece}.]
One bounds,
remembering that the variables $\varphi_k$ are martingale differences
and using
the domination \eqref{eqnphnph},
\begin{eqnarray*}
\bigl|\mathbb{E}(\varphi_{n,k}\vert \mathcal{B}_{k-1})\bigr|&=& \bigl|
\mathbb{E}(\varphi_k-\varphi_{n,k}\vert
\mathcal{B}_{k-1})\bigr|
\\
&\leq& \int_{n^\gamma}^\infty\mathbb{P}\bigl(|
\varphi_k|>t\vert \mathcal {B}_{k-1}\bigr)\,\mathrm{d}t +
n^\gamma\mathbb{P}\bigl(|\varphi_k|> n^\gamma\vert
\mathcal{B}_{k-1}\bigr)
\\
&\leq & \int_{n^\gamma}^\infty\mathbb{P}(\varphi>t)\,
\mathrm{d}t + n^\gamma\mathbb{P}\bigl(\varphi> n^\gamma\bigr) =
\mathbb{E}(\varphi\,\mathbf{1}_{\{\varphi>n^\gamma\}}),
\end{eqnarray*}
and this right-hand side converges to $0$ when $n$ goes to infinity since
the dominating function $\varphi$ is integrable. One deduces the bounds
\[
\frac{1}n |T_n-\overline T_n|\leq \mathbb{E}(
\varphi \mathbf{1}_{\{\varphi>n^\gamma\}}),
\]
with a right-hand side also converging to $0$.
Hence, one can find an integer $n_0=n_0(p,\varepsilon,\varphi)$ such that,
for $n\geq n_0$,
the event $A_{3,n}$ is empty.
We just set $C_{3,n}=0$ when $n\geq n_0$ and $C_{3,n}=1$ otherwise.

\item[\textit{Fourth piece}.] We set
$\overline Q_n:=\sum_{1\leq k\leq n}\overline\varphi_{n,k}^2$,
$p_0:=\min(p,2)$, and $M\geq1$ to be the smallest integer
such that $M\geq\frac{p}{2(1-\gamma)}$.
According to the Burkholder inequality (see \cite{HH80}),
since $(\overline\varphi_{n,k})_{1\leq k\leq n}$ is a martingale difference,
there exists
a constant $D_M$, which depends only on $M$, such that
\[
D_M^{-1} \mathbb{E}\bigl(\overline{Q}_n^M
\bigr) \leq\mathbb{E}\bigl(\overline{T}_n^{2M}\bigr) \leq
D_M \mathbb{E}\bigl(\overline{Q}_n^M\bigr).
\]
One computes then, using Chebyshev's inequality,
%
\begin{equation}
\label{eqnpanetneqn} \mathbb{P}(A_{4,n})\leq n^{-2M} \mathbb{E}
\bigl(\overline{T}_n^{2M}\bigr)\leq D_M
n^{-2M} \mathbb{E}\bigl(\overline{Q}_n^M\bigr).
\end{equation}
We expand now $\mathbb{E}(\overline{Q}_n^{M})$ as a sum of terms of the form
$\mathbb{E}(\overline{\varphi}_{n,k_1}^{2q_1}\cdots\overline{\varphi}
_{n,k_\ell}^{2q_\ell})$
with $1\leq\ell\leq M$,
$q_1,\ldots, q_{\ell}\geq1$, $q_1 +\cdots+q_{\ell}=M$ and $1\leq
k_1<\cdots<k_{\ell}\leq n$.
Using the bounds, for $1\leq k\leq n$ and $q\geq1$,
\[
\overline{\varphi}_{n,k}^{2q}\leq\bigl(2n^\gamma
\bigr)^{2q-p_0}| \overline{\varphi}_{n,k}|^{p_0},
\]
and, using the domination \eqref{eqnphnph}, one bounds each term in
the sum
\[
\mathbb{E}\bigl(\overline\varphi_{n,k_1}^{2q_1}\cdots\overline
\varphi _{n,k_\ell}^{2q_\ell}\bigr) \leq 4^{M}n^{2M\gamma-\ell p_0\gamma}
\mathbb{E} \bigl(\varphi^{p_0}\bigr)^\ell.
\]
For each value of $\ell\leq M$,
the number of such terms is bounded by
$M^{\ell} n^{\ell} $. Summing all these bounds, one gets, since
$\gamma
p_0>\min(\frac{p+1}{2},\frac{p+1}{p})>1$,
\begin{eqnarray*}
\mathbb{E}\bigl(\overline{Q}_n^M\bigr) &\leq & \sum
_{1\leq\ell\leq M} (4M)^M \mathbb{E} \bigl(
\varphi^{p_0}\bigr)^\ell n^{2M\gamma-\ell p_0\gamma
+\ell}
\\
&\leq & c_{p,\varphi} n^{2M\gamma},
\end{eqnarray*}
where $c_{p,\varphi}=4^M M^{M+1}\max(1,\mathbb{E} (\varphi^{p_0})^M)$.
Plugging this inside \eqref{eqnpanetneqn}, one gets
\[
\mathbb{P}(A_{4,n})\leq C_{4,n}:= c_{p,\varphi}
D_{M } n^{-2(1-\gamma)M}
\]
and
\[
\sum_{n\geq1}n^{p-2}C_{4,n}=c_{p,\varphi}
D_{M } \sum_{n\geq1}n^{p-2-2(1-\gamma)M },
\]
which is finite since $M\geq\frac{p}{2(1-\gamma)}$.\quad\qed
\end{longlist}
\noqed\end{pf*}

\begin{Rmq}
As we have seen in this proof, assumption \eqref{eqnphnph} in
Theorem~\ref{thhsurobmar} implies that there exists a constant
$C:=\mathbb{E}|\varphi|^p$
such that, for all $n\geq1$,
%
\begin{equation}
\label{eqnphnpbn} \mathbb{E}\bigl(|\varphi_n|^p\vert
\mathcal{B}_{n-1}\bigr) \leq C.
\end{equation}
However, the conclusion of Theorem~\ref{thhsurobmar} is no more true
if we replace assumption \eqref{eqnphnph} by \eqref{eqnphnpbn}.
Here is a counterexample. Choose $\varphi_n$ to be symmetric independent
random variables such that, for $3^{i-1}< n\leq3^i$,
$\varphi_n$ takes values in the set $\{ -3^i, 0, 3^i\}$ and
$\mathbb{P} (\varphi_n=\pm3^i)=3^{-pi}$. For these variables,
the conclusion of Theorem~\ref{thhsurobmar} does not hold. This is
essentially due to the fact
that the series $\sum_{n\geq1}n^{p-2}\mathbb{P}(\exists k\leq n\vert
  |\varphi
_k|\geq n )$ diverge
(the details are left to the reader since we will not use this example).
\end{Rmq}

When the martingale difference is uniformly bounded,
one has a much better large deviation estimate
than \eqref{eqnpsncn} due to Azuma in \cite{Az67}.

\begin{Fact}[(Azuma)]
\label{thazuma}
Let $(\varphi_n)_{n\geq1}$ be a martingale difference
and $S_n:=\varphi_1+\cdots+\varphi_n$ the corresponding martingale.
If $|\varphi_n|\leq a<\infty$ for all $n\geq1$, then one has
for all $n\geq1$, $\varepsilon>0$,
%
\begin{equation}
\label{eqnpsnen} \mathbb{P}(S_n\geq n\varepsilon)\leq
e^{-({n\varepsilon^2})/({2a^2})}.
\end{equation}
\end{Fact}

\begin{pf}
We recall Azuma's proof since it is very short. Assume $a=1$.
Using the convexity of the exponential function, one bounds, for all
$x$ in $[-1,1]$,
$e^{\varepsilon x}\leq\cosh(\varepsilon)+ x\sinh(\varepsilon)\leq
e^{{\varepsilon^2}/{2}}+x\sinh(\varepsilon)$. Hence,
for all $k\geq1$, one has $\mathbb{E}(e^{\varepsilon\varphi_k}\vert
\mathcal{B}_{k-1})\leq
e^{{\varepsilon^2}/{2}}$,
and, by Chebyshev's inequality,
\[
\mathbb{P}(S_n\geq n\varepsilon)\leq e^{-n\varepsilon^2} \mathbb {E}
\bigl(e^{\varepsilon S_n}\bigr) \leq e^{-n\varepsilon^2}\bigl(e^{{\varepsilon^2}/{2}}
\bigr)^n =e^{-({n\varepsilon^2})/{2}}.
\]
\upqed\end{pf}

\subsection{Central limit theorem for martingales}
\label{seccltmar}
%
In this section, we briefly recall the martingale central limit
theorem, which is due to Brown.

Let $(\Omega,\mathcal{B},\mathbb{P})$ be a probability space,
$(p_n)_{n\geq1}$
be a sequence of positive integers and,
for $n\geq1$, let
\[
\mathcal{B}_{n,0}\subset\cdots\subset\mathcal{B}_{n,p_n}
\]
be sub-$\sigma$-algebras of $\mathcal{B}$.

Let $E$ be a finite-dimensional normed real vector space.
We want to define the Gaussian laws $N_\Phi$ on $E$.
Such a law is completely determined by its \textit{covariance $2$-tensor}
$\Phi$.
If we fix a Euclidean structure on $E$, this \textit{covariance
$2$-tensor} is nothing but
the \textit{covariance matrix} of $N_\Phi$. Here are the precise definitions.

We denote by $S^2E$ the space of symmetric $2$-tensors of $E$.
Equivalently, $S^2E$ is the space
of quadratic forms on the dual space $E^*$.
The \textit{linear span} of a symmetric $2$-tensor $\Phi$
is the smallest vector subspace $E_\Phi\subset E$
such that $\Phi$ belongs to $\mathrm{S}^2E_\Phi$.
A $2$-tensor $\Psi\in\mathrm{S}^2E$ is \textit{nonnegative}
(which we write $\Psi\geq0$) if it is nonnegative as a quadratic form on
the dual space $E^*$.
For every $v$ in $E$, we set $v^2:=v\otimes v\in\mathrm{S}^2E$,
and we denote by
\[
B_\Phi:=\bigl\{v\in E_\Phi\vert \Phi - v^2
\mbox{ is nonnegative}\bigr\}
\]
the \textit{unit ball} of $\Phi$.
For any nonnegative symmetric $2$-tensor $\Phi\in\mathrm{S}^2E$, we let
$N_\Phi$ be the centered Gaussian law on $E$
with covariance $2$-tensor $\Phi$, that is, such that
\[
\Phi=\int_E v^2 \,\mathrm{d}N_\Phi(v).
\]
For instance, $N_\Phi$ is a Dirac mass at $0$ if
and only if
$\Phi=0$ if and only if $E_\Phi=\{ 0\}$.

The following theorem is due to Brown in
\cite{Bro71} (see also \cite{HH80}).

\begin{Fact}[(Brown martingale central limit theorem)]
\label{thcltmar}
\textup{For $1\leq k\leq p_n$, let $\varphi_{n,k}\dvtx \Omega\rightarrow E$
be square-integrable random variables such that
%
\begin{equation}
\label{eqnenkphnk} \mathbb{E}(\varphi_{n,k}\vert \mathcal{B}_{n,k-1})=0.
\end{equation}
We assume
that the $\mathrm{S}^2E$-valued random variables
%
\begin{equation}
\label{eqnwntdph} W_n:=\sum_{1\leq k\leq p_n}
\mathbb{E}\bigl(\varphi_{n,k}^2\vert \mathcal
{B}_{n,k-1}\bigr) \qquad \mbox{converge to $\Phi$ in probability},
\end{equation}
and that, for all $\varepsilon>0$,
%
\begin{equation}
\label{eqnwnepstd0} W_{\varepsilon,n} := \sum_{1\leq k\leq p_n}
\mathbb{E}\bigl(\varphi_{n,k}^2\mathbf{1}_{\{\|\varphi_{n,k}\|\geq
\varepsilon\}}
\vert \mathcal{B}_{n,k-1}\bigr) \mathop{\longrightarrow}_{n\rightarrow\infty}0
\qquad \mbox{in probability}.\hspace*{-8pt}
\end{equation}
Then the sequence
$S_n:=\sum_{1\leq k\leq p_n}\varphi_{n,k}$ converges in law
toward $N_\Phi$.}
\end{Fact}

Under the same assumptions, the sequence $S_n$ also satisfies a law of
the iterated logarithm,
that is, almost surely, the set of cluster points of the sequence
$\frac{S_n}{\sqrt{2\Phi  n\log\log n}}$ is equal to the unit ball
$B_\Phi$
(indeed the sequence $S_n$ satisfies
an invariance principle; see \cite{HH80}, Chapter~4).

Assumption \eqref{eqnwnepstd0} is called Lindeberg's condition.

We recall that a sequence $X_n$ of random variables
converges to $X_\infty$ in probability, if, for all $\varepsilon>0$,
$\mathbb{P}(|X_n-X_\infty|\geq\varepsilon)\displaystyle{\mathop{\longrightarrow}_{n\rightarrow\infty}} 0$.

\section{Limit theorems for cocycles}
\label{seccencoc}
%
In this section, we state various limit theorems for
cocycles and we explain how to deduce them
from the limit theorems for martingales that we discussed in
Chapter~\ref{seclimmar}.

\subsection{Complete convergence for functions}
\label{secllnfon}
%
In this section, we prove a large deviations estimate
in the law of large numbers for functions
over Markov--Feller chains.

Let $X$ be a compact metrizable space and $\mathcal{C}^0(X)$ be the
Banach space
of continuous functions on $X$.
Let $P\dvtx \mathcal{C}^0(X)\rightarrow\mathcal{C}^0(X)$ be a
Markov--Feller operator,
that is, a bounded operator such that $\| P\|\leq1$, $P1=1$ and such that
$Pf\geq0$ for all functions $f\geq0$.
Such a Markov--Feller operator can be seen alternatively as a
weak-$*$ continuous map $x\mapsto P_x$ from $X$ to the set of
probability measures on $X$,
where $P_x$ is defined\vspace*{1pt} by $P_x(f)=(Pf)(x)$ for all $f$ in $\mathcal{C}^0(X)$.
We denote by $\underline X$ the compact set $\underline X=X^{\mathbb{N}}$
of infinite sequences $\underline x=(x_0,x_1,x_2,\ldots)$.
For $x$ in $X$, we denote by $\mathbb{P}_x$ the Markov probability
measure on
$\underline{X}$,
that is, the law of the trajectories of the Markov chain starting from
$x$ associated to $P$.

Given a continuous function $\varphi$ on $X$,
we define its \textit{upper average} by
\[
\ell^+_\varphi={\sup_{\nu}} \int_G
\varphi(x)\,\mathrm{d}\nu(x)
\]
and \textit{lower average} by
\[
\ell^-_\varphi:= \inf_{\nu} \int_G
\varphi(x)\,\mathrm{d}\nu(x),
\]
where the supremum and the infimum are taken over all the $P$-invariant
probability measures $\nu$ on $X$.
We say $\varphi$ \textit{has unique average} if $\ell^+_\varphi=\ell
^-_\varphi$.

According to the Breiman law of large numbers
in \cite{Brei60} (see also \cite{BQRW}),
for such a $\varphi$, for any $x$ in $X$,
for $\mathbb{P}_x$-almost every $\underline{x}$ in $\underline{X}$, the sequence
$
\frac{1}{n}\sum_{k=1}^{n}\varphi(x_k)$ converges to $\ell
^+_\varphi=\ell^-_\varphi$.
The following proposition is a large deviations estimate for the
Breiman law of large numbers.

\begin{Prop}
\label{prolardevbre}
Let $X$ be a compact metrizable space, and
$P$ be a Markov--Feller operator on $X$.
Let $\varphi$ be a continuous function on $X$ with upper average $\ell
^+_\varphi$
and lower average $\ell^-_\varphi$.
Then, for all $\varepsilon>0$, there exist constants $A>0$, $\alpha>0$
such that
%
\begin{equation}
\label{eqnlardevbre} \mathbb{P}_x \Biggl(\Biggl\{ \underline x\in
\underline X\Big\vert \frac{1}n\sum_{k=1}^n
\varphi(x_k)\notin\bigl[\ell^-_\varphi - \varepsilon,
\ell^+_\varphi + \varepsilon\bigr]\Biggr\} \Biggr) \leq A
e^{-\alpha n},
\end{equation}
for all $n\geq1$ and all $x$ in $X$.
\end{Prop}

Note that $\ell^-_\varphi=\ell^+_\varphi$ as soon as $P$ is
uniquely ergodic,
that is, as soon as there exists only one $P$-invariant Borel
probability measure $\nu$ on $X$.

\begin{pf*}{Proof of Proposition~\protect\ref{prolardevbre}}
We assume $\|\varphi\|_\infty=\frac{1}2$.
We introduce, for $1\leq\ell\leq n$, the bounded functions $\Psi_n$ and
$\Psi_{\ell,n}$ on $\underline X$
given, for $\underline x$ in $ \underline X$, by
\[
\Psi_n(\underline x)=\varphi(x_n) \quad\mbox{and}\quad
\Psi_{\ell,n}(\underline x)=\bigl(P^{\ell}\varphi\bigr)
(x_{n-\ell}),
\]
so that, for $x$ in $X$,
\[
\Psi_{\ell,n}=\mathbb{E}_x(\Psi_n\vert
\mathcal{X}_{n-\ell}) \qquad \mathbb{P}_x\mbox{-a.s.},
\]
where $\mathcal{X}_{n}$ is the $\sigma$-algebra on $\underline X$ spanned
by the functions $\underline x\mapsto x_k$ with $k\leq n$.
On one hand,
one has the uniform convergence
%
\begin{equation}
\label{eqnmsumjmph} \max\Biggl(\ell^+_\varphi,\frac{1}{m}\sum
_{j=1}^m P^{j}\varphi \Biggr)\mathop{
\longrightarrow}_{m\rightarrow\infty}\ell ^+_\varphi
\end{equation}
in $\mathcal{C}^0(X)$.
Hence, we can fix $m$ such that,
for all $x\in X$,
\[
\label{eqnsumpjph} \frac{1}{ m}\sum_{j=1}^mP^j
\varphi(x) \leq\ell^+_\varphi+ \frac{\varepsilon}{4}.
\]
Then,
for all $n\geq1$ and $\underline x\in\underline X$, one has
\[
\label{eqnsumsumpsi1} \frac{1}{n  m}\sum_{k=m+1}^{m+n}
\sum_{j=1}^m\Psi _{j,k+j}(
\underline x) \leq\ell^+_\varphi+ \frac{\varepsilon}{4}.
\]
In particular, if $n\geq n_0:=\frac{4m}{\varepsilon}$, one also has
%
\begin{equation}
\label{eqnsumsumpsi2} \frac{1}{n  m}\sum_{k=m+1}^{m+n}
\sum_{j=1}^m\Psi_{j,k}(
\underline x)\leq \ell^+_\varphi+ \frac{\varepsilon}{2}.
\end{equation}

On the other hand,
for all $1\leq j\leq m$, $x\in X$, by Azuma's bound \eqref{eqnpsnen}
and the equalities, for $k\geq j$, $\mathbb{E}_x(\Psi_{j-1,k}-\Psi
_{j,k}\vert
\mathcal{X}_{k-j})=0$, one has
\[
\mathbb{P}_x\Biggl(\Biggl\{\underline x\in\underline X\Big\vert \Biggl|
\frac{1}{n} \sum_{k=m+1}^{m+n}\bigl(
\Psi_{j-1,k}(\underline x) -\Psi_{j,k}(\underline x)\bigr)\Biggr| \geq
\frac{\varepsilon}{4m}\Biggr\}\Biggr)\leq e^{-({n\varepsilon^2})/({32m^2})}.
\]
Adding these bounds, one gets, for all $1\leq j\leq m$, $x\in X$,
\[
\mathbb{P}_x\Biggl(\Biggl\{\underline x\in\underline X\Big\vert \Biggl|
\frac{1}{n} \sum_{k=m+1}^{m+n}\bigl(
\Psi_k(\underline x) -\Psi_{j,k}(\underline x)\bigr)\Biggr| \geq
\frac{\varepsilon}{4}\Biggr\}\Biggr)\leq m e^{-({n\varepsilon^2})/({32m^2})},
\]
and hence
\[
\label{eqnbepsikpsijk} \mathbb{P}_x\Biggl(\Biggl\{\underline{x}\in
\underline{X}\Big\vert \Biggl|\frac{1}{n}\sum_{k=m+1}^{m+n}
\Biggl(\Psi_k(\underline{x}) - \frac{1}{m}\sum
_{j=1}^m \Psi_{j,k}(\underline{x})
\Biggr)\Biggr| \geq\frac{\varepsilon}{4}\Biggr\}\Biggr)\leq m^2
e^{-({n\varepsilon^2})/({32m^2})}.
\]
Combining this formula with \eqref{eqnsumsumpsi2},
one gets the desired bound,
\[
\mathbb{P}_x\Biggl(\Biggl\{\underline{x}\in\underline{X}\Big\vert
\frac{1}{n}\sum_{k=1}^{n}
\Psi_k(\underline{x})\geq \ell^+_\varphi+\varepsilon\Biggr\}
\Biggr)\leq m^2 e^{-({n\varepsilon^2})/({32m^2})},
\]
for all $n\geq n_0$ and $x\in X$.
\end{pf*}

\subsection{Complete convergence for cocycles}
\label{seccllncoc}
%
In this section, we prove the complete convergence
in the law of large numbers for cocycles
over $G$-spaces.

Let
$G$ be a second countable locally compact group
acting continuously on a compact second countable topological space $X$.
Let $\mu$ be a Borel probability measure on $G$.

We denote by $(B,\mathcal{B},\beta)$
the associated one-sided
Bernoulli space,
that is, $B=G^{\mathbb N^*}$
is the set of sequences $b=(b_1,\ldots,b_n,\ldots)$ with $b_n$ in $G$,
$\mathcal B$ is the product
$\sigma$-algebra of the Borel $\sigma$-algebras of $G$, and $\beta$ is
the product measure
$\mu^{\otimes\mathbb N^*}$.
For $n\geq1$, we denote by $\mathcal{B}_n$
the $\sigma$-algebra spanned by the $n$ first coordinates $b_1,\ldots,b_n$.

We will apply the results of Section~\ref{secllnfon}
to the \textit{averaging operator},
that is, the
Markov--Feller operator
$P=P_\mu\dvtx \mathcal{C}^0(X)\rightarrow\mathcal{C}^0(X)$ whose
transition probabilities are
given by
$P_x=\mu*\delta_x$ for all $x$ in $X$.
For every $x$ in $X$, the Markov measure $\mathbb{P}_x$
is the image of $\beta$ by the map
\[
B\rightarrow\underline{X};\qquad b\mapsto(x, b_1 x,
b_2b_1x,b_3b_2b_1x,
\ldots).
\]
We denote by $\mu^{*n}$ the $n$th-convolution power $\mu
*\cdots
*\mu$.

Let $E$ be a finite-dimensional normed real vector space
and $\sigma$ a continuous function
$\sigma\dvtx G\times X\rightarrow E$.
This function $\sigma$ is said to be a \textit{cocycle}
if one has
%
\begin{equation}
\label{eqncocsighx} \sigma\bigl(gg',x\bigr)=\sigma\bigl(g,g'x
\bigr) + \sigma\bigl(g',x\bigr)\qquad \mbox{for any
$g,g'\in G$, $x\in X$.}
\end{equation}

We introduce the
\textit{sup-norm} function $\sigma_{\sup}$.
\index{sup-norm}
It is given, for $g$ in $G$, by
%
\begin{equation}
\label{eqnkagsigx} \sigma_{\sup}(g)= \sup_{x\in X} \bigl\|
\sigma(g,x)\bigr\|.
\end{equation}
We assume that this function $\sigma_{\sup}$ is integrable
%
\begin{equation}
\label{eqnintsisup} \int_G\sigma_{\sup}(g) \,
\mathrm{d}\mu(g)<\infty.
\end{equation}
Recall a Borel probability measure $\nu$ on $X$ is said to be $\mu$-stationary if $\mu*\nu=\nu$, that is, if it is $P_\mu$-invariant.
When $E=\mathbb{R}$, we define the \textit{upper average} of $\sigma$ by
\[
\sigma^+_\mu=\sup_{\nu} \int_{G\times X}
\sigma(g,x)\,\mathrm{d}\mu(g)\,\mathrm{d}\nu(x),
\]
and the \textit{lower average}
\[
\sigma^-_\mu=\inf_{\nu} \int_{G\times X}
\sigma(g,x)\,\mathrm {d}\mu(g)\,\mathrm{d}\nu(x),
\]
where the supremum and the infimum are taken over all the $\mu$-stationary
probability measures $\nu$ on $X$.
We say that $\sigma$ \textit{has unique average} if
the averages
do not depend on the choice of
the $\mu$-stationary
probability measure $\nu$, that is, if $\sigma_\mu^+=\sigma_\mu^-$.
In this case, these functions satisfy also a law of large numbers,
that is, under assumption
\eqref{eqnintsisup} if $\sigma$ has unique average,
for any $x$ in $X$, for $\beta$-almost every $b$ in $B$, the sequence
$\sum_{k=1}^n\frac{\sigma(b_k,b_{k  -  1}\cdots b_1x)}{n}$
converges to $\sigma_\mu$ (see \cite{BQRW}, Chapter~2).

Proposition \ref{prollncoc}
is an analog of the Baum--Katz theorem for these functions.
For $p=2$, it says that, when $\sigma_{\sup}$ is square integrable, this
sequence converges completely.

\begin{Prop}
\label{prollncoc}
Let $G$ be a locally compact group, $X$ a compact metrizable $G$-space,
$\mu$ a Borel probability measure
on $G$ and $p>1$.
Let $\sigma \dvtx  G\times X\rightarrow\mathbb{R}$ be a continuous function
such that
$\sigma_{\sup}$ is $\mathrm{L}^p$-integrable.
Let $\sigma^+_\mu$ and $\sigma^-_\mu$ be its
upper and lower average.
Then, for any $\varepsilon>0$, there exist constants $D_n$ such that
\[
\sum_{n\geq1}n^{p-2}D_n<\infty,
\]
and, for $n\geq1$, $x\in X$,
\[
\label{eqnbesumsigkg1x} \beta\Biggl(\Biggl\{b\in B\Big\vert \sum
_{k=1}^n\frac{\sigma(b_k,b_{k  -
1}\cdots b_1x)}{n} \notin\bigl[
\sigma^-_\mu-\varepsilon,\sigma^+_\mu+\varepsilon\bigr]
\Biggr\} \Biggr)\leq D_n.
\]
In particular, when $\sigma$ is a cocycle, one has, for $n\geq1$,
$x\in X$,
%
\begin{equation}
\label{eqnbesigng1x} \mu^{*n}\biggl(\biggl\{g\in G\Big\vert \frac{\sigma(g,x)}{n}
\notin\bigl[\sigma^-_\mu-\varepsilon,\sigma^+_\mu+
\varepsilon\bigr]\biggr\} \biggr)\leq D_n.
\end{equation}
\end{Prop}

The fact that the constants $D_n$ do not depend on $x$
will be important for our applications.

\begin{pf*}{Proof of Proposition~\protect\ref{prollncoc}}
According to Proposition~\ref{prolardevbre}, the conclusion
of Proposition~\ref{prollncoc} is true when the function $\sigma$
does not depend on the variable~$g$.
Hence, it is enough to prove Proposition~\ref{prollncoc}
for the continuous function
$\sigma'$ on $G\times X$
given, for $g$ in $G$ and $x$ in $X$, by
\[
\sigma'(g,x)=\sigma(g,x)-\int_G\sigma(g,x)
\,\mathrm{d}\mu(g).
\]
By construction, the sequence of functions $\varphi_n$ on $B$ given, for
$b$ in $B$, by
\[
\varphi_n(b)=\sigma'(b_n,b_{n  -  1}
\cdots b_1x)
\]
is a martingale difference.
Hence, our claim follows from Theorem~\ref{thhsurobmar}
since the functions $\varphi_n$ satisfy the domination \eqref{eqnphnph}:
for $n\geq1$, $t>0$,
\[
\label{eqnphnph2} \mathbb{E}(\mathbf{1}_{\{|\varphi_n|>t\}}\vert \mathcal{B}_{n-1})
\leq \mu\bigl(\bigl\{g\in G\vert \sigma_{\sup}(g)+M>t\bigr\}\bigr),
\]
where $M$ is the constant $M:=\int_G\sigma_{\sup}(g)\,\mathrm{d}\mu(g)$.
\end{pf*}

\subsection{Central limit theorem for centerable cocycles}
\label{seccltcoc}
%
In this section, we explain how to deduce the central limit theorem for
centerable cocycles from the central limit theorem for martingales.

Let $\sigma \dvtx  G\times X\rightarrow E$ be a continuous cocycle.
When the function $\sigma_{\sup}$ is $\mu$-integrable, one defines the
\textit{drift} or \textit{expected increase} of $\sigma$: it is the
continuous function
$X\rightarrow E; x\mapsto\int_G\sigma(g,x)\,\mathrm{d}\mu(g)$. One
says that $\sigma$ has
\index{constant drift} \textit{constant drift} if the drift is a constant
function:
%
\begin{equation}
\label{eqnintgsigxmug} \int_G\sigma(g,x)\,\mathrm{d}\mu(g)=
\sigma_\mu.
\end{equation}
One says that $\sigma$ is
\textit{centered} if the drift is a null function.

A continuous cocycle $\sigma \dvtx  G\times X\rightarrow E$ is said
to be \textit{centerable} if it is the sum
%
\begin{equation}
\label{eqnspecoc} \sigma(g,x)=\sigma_0(g,x) +\psi(x)-\psi(gx)
\end{equation}
of a cocycle $\sigma_0(g,x) $ with constant drift $\sigma_\mu$ and
of a
coboundary $\psi(x)-\psi(gx)$ given by a continuous function $\psi
\in
\mathcal{C}^0(X)$.
A centerable cocycle always has a unique average:
for any $\mu$-stationary probability $\nu$ on $X$,
one has
\[
\int_{G\times X} \sigma(g,x)\,\mathrm{d}\mu(g)\,\mathrm{d}\nu (x)=
\sigma_\mu.
\]

Here is a trick to reduce the study of a cocycle with constant drift
$\sigma_\mu$
to one which is centered. Replace $G$ by $G':=G\times\mathbb{Z}$
where $\mathbb{Z}$ acts trivially on $X$,
replace $\mu$ by $\mu':=\mu\otimes\delta_1$, so that any $\mu$-stationary
probability
measure on $X$ is also
$\mu'$-stationary, and replace $\sigma$ by the cocycle
%
\begin{equation}
\label{eqncenteredcocycle} \sigma' \dvtx G'\times X\rightarrow E
\qquad \mbox{given by } \sigma'\bigl((g,n),x\bigr)=\sigma(g,x)-n
\sigma_\mu.
\end{equation}

A centerable cocycle $\sigma$ is said to \textit{have unique
covariance}
$\Phi_\mu$
if
%
\begin{eqnarray}
\Phi_\mu &:= & \int_{G\times X}\bigl(
\sigma_0(g,x)-\sigma_\mu\bigr)^2\,\mathrm {d}
\mu(g)\,\mathrm{d}\nu(x)
\nonumber
\\[-8pt]
\label{eqnphmusisimu}
\\[-8pt]
\eqntext{\mbox{does not depend on the choice of the $\mu$-stationary probability measure $
\nu$},}
\end{eqnarray}
where $\sigma_0$ is as in \eqref{eqnspecoc}.
This covariance $2$-tensor $\Phi_\mu\in\mathrm{S}^2E$ is nonnegative.

\begin{Rmq}
This assumption does not depend on the choice of $\sigma
_0$. More precisely, if $\sigma_0$ and $\sigma_1$ are cohomologous
centered cocycles,
for any $\mu$-stationary Borel probability measure $\nu$ on $X$, one has
%
\begin{equation}
\label{eqnsi01} \int_{G\times X}\sigma_0(g,x)^2
\,\mathrm{d}\mu(g)\,\mathrm{d}\nu (x)=\int_{G\times X}\sigma
_1(g,x)^2\,\mathrm{d}\mu(g)\,\mathrm{d}\nu(x).
\end{equation}
Indeed, since $\sigma_0$ and $\sigma_1$ are centered and cohomologous,
we may write, for any $g,x$,
$\sigma_1(g,x)=\sigma_0(g,x)+\psi(x)-\psi(gx)$ where $\psi$ is a
continuous function on $X$ and $P_\mu\psi=\psi$.
Now, the difference between the two sides of \eqref{eqnsi01} reads as
%
\begin{equation}
\label{eqnsipsi} 2\int_{G\times X}\sigma_0(g,x)\psi(gx)
\,\mathrm{d}\mu(g)\,\mathrm {d}\nu(x).
\end{equation}
By ergodic decomposition, to prove this is $0$, one can assume $\nu$ is
$\mu$-ergodic. In this case, since $P_\mu\psi=\psi$,
$\psi$ is constant $\nu$-almost everywhere and \eqref{eqnsipsi} is
proportional to
$ \int_{G\times X}\sigma_0(g,x)\,\mathrm{d}\mu(g)\,
\mathrm{d}\nu(x)$,
which is $0$ by assumption.
\end{Rmq}

\begin{Thm}[(Central limit theorem for centerable cocycles)]
\label{thcltcoc}
Let $G$ be a locally compact group, $X$ a compact metrizable $G$-space,
$E$ a finite-dimensional real vector space, and $\mu$ a Borel
probability measure
on $G$.
Let $\sigma \dvtx  G\times X\rightarrow E$ be a continuous cocycle such that
$\int_G\sigma_{\sup}(g)^2\,\mathrm{d}\mu(g)<\infty$.
Assume that $\sigma$ is centerable with average $\sigma_\mu$
and has a unique covariance $\Phi_\mu$, that is, $\sigma$ satisfies
\eqref{eqnspecoc} and \eqref{eqnphmusisimu}.
Let $N_\mu$ be
the Gaussian law on $E$ whose covariance $2$-tensor is~$\Phi_\mu$.

Then, for any bounded continuous
function $\psi$ on $E$, uniformly for $x$ in $X$,
one has
%
\begin{equation}
\label{eqntclcoc} \int_{G}\psi \biggl(\frac{\sigma(g,x)-n\sigma_\mu}{
\sqrt{n}}
\biggr)\,\mathrm{d}\mu^{*n}(g) \mathop{\longrightarrow}_{n\rightarrow\infty}
\int_{ E}\psi(v)\,\mathrm{d}N_\mu(v).
\end{equation}
\end{Thm}

Note that hypothesis \eqref{eqnphmusisimu} is automatically satisfied
when there exists a unique $\mu$-stationary Borel probability measure
$\nu$ on $X$.

\begin{Rmqs}
When $E=\mathbb{R}^d$, the covariance $2$-tensor $\Phi_\mu$ is nothing
but the \textit{covariance matrix} of the random variable $\sigma_0$
on $(G\times X,\mu\otimes\nu)$.

The conclusion in Theorem~\ref{thcltcoc} is not correct
if one does not assume the cocycle $\sigma$ to be centerable.
\end{Rmqs}

\begin{pf*}{Proof of Theorem~\protect\ref{thcltcoc}}
We will deduce Theorem~\ref{thcltcoc} from
the central limit Theorem \ref{thcltmar} for martingales.

As in the previous sections,
let $(B,\mathcal{B},\beta)$ be the Bernoulli space with alphabet
$(G,\mu)$.
We want to prove that, for any sequence $x_n$ on $X$,
the laws of the random variables $S_n$ on $B$
given, for $b$ in $B$, by
\[
S_n(b):=\frac{1}{\sqrt{n}}\bigl(\sigma(b_n\cdots
b_1, x_n)-n\sigma_\mu\bigr)
\]
converge to $N_\mu$.

Since the cocycle $\sigma$ is centerable, one can write
$\sigma$ as the sum of two cocycles
$\sigma=\sigma_0+\sigma_1$ where $\sigma_0$ has constant drift and
where $\sigma_1$ is a coboundary. In particular, the cocycle
$\sigma_1$ is uniformly bounded and does not play any role in the
limit~\eqref{eqntclcoc}.
Hence, we can assume $\sigma=\sigma_0$.
Using the trick \eqref{eqncenteredcocycle},
we can assume that $\sigma_\mu=0$, that is, $\sigma$ is a centered cocycle.

We want to apply the martingale central limit Theorem
\ref{thcltmar}
to the sub-$\sigma$-algebras $\mathcal{B}_{n,k}=\mathcal{B}_k$
spanned by $b_1,\ldots,b_k$ and to
the triangular array of random variables
$\varphi_{n,k}$ on $B$ given by, for $b$ in $B$,
\[
\varphi_{n,k}(b)=\frac{1}{\sqrt{n}}\sigma(b_k,b_{k-1}
\cdots b_1 x_n) \qquad \mbox{for } 1\leq k\leq n.
\]
Since, by the cocycle property \eqref{eqncocsighx}, one
has
\[
S_n=\sum_{1\leq k\leq n}\varphi_{n,k},
\]
we just have to check
that the three assumptions of Theorem~\ref{thcltmar}
are satisfied with $\Phi=\Phi_\mu$.
We keep the notation
$W_n$ and $W_{\varepsilon,n}$ of this theorem.

First, since the function $\sigma_{\sup}$ is square integrable, the
functions
$\varphi_{n,k}$ belong to $\mathrm{L}^2(B,\beta)$, and,
by equation \eqref{eqnintgsigxmug}, assumption
\eqref{eqnenkphnk} is satisfied: for $\beta$-almost all $b$ in $B$,
\[
\mathbb{E}(\varphi_{n,k}\vert \mathcal{B}_{k-1})= \int
_G\sigma(g,b_{k-1}\cdots b_1x_n)
\,\mathrm{d}\mu(g)=0.
\]

Second, we introduce the continuous function on $X$,
\[
x\mapsto M(x)=\int_G \sigma(g,x)^2\,
\mathrm{d}\mu(g)
\]
and we compute, for $\beta$-almost all $b$ in $B$,
\[
W_{n}(b)=\frac{1}{n}\sum_{1\leq k\leq n}M(b_{k-1}\cdots b_1x_n).
\]
According to Proposition~\ref{prolardevbre},
since $\sigma$ has a unique covariance $\Phi_\mu$,
the sequence $W_n$ converges to $\Phi_\mu$ in probability,
that is, assumption
\eqref{eqnwntdph} is satisfied.

Third, we introduce, for $\lambda>0$,
the continuous function on $X$
\[
x\mapsto M_\lambda(x)=\int_G
\sigma(g,x)^2\mathbf{1}_{\{\|\sigma(g,x)\|
\geq\lambda\}}\,\mathrm{d} \mu(g)
\]
and the integral
\[
I_\lambda:=\int_G\sigma_{\sup}^2(g)
\mathbf{1}_{\{\sigma_{\sup}(g)\geq\lambda\}} \,\mathrm{d}\mu(g),
\]
we notice that
\[
M_\lambda(x)\leq I_\lambda\mathop{\longrightarrow}_{\lambda\rightarrow\infty}0,
\]
and we compute, for $\varepsilon>0$ and $\beta$-almost all $b$ in $B$,
\[
W_{\varepsilon,n}(b)=\frac{1}{n}\sum_{1\leq k\leq n}M_{\varepsilon
\sqrt{n}}(b_{k-1}\cdots b_1x_n)
\leq I_{\varepsilon\sqrt{n}}\mathop{\longrightarrow}_{n\rightarrow\infty}0.
\]
In particular, the sequence $W_{\varepsilon,n}$ converges to $0$ in
probability, that is, Lindeberg's condition
\eqref{eqnwnepstd0} is satisfied.

Hence, by Fact~\ref{thcltmar}, the laws of $S_n$ converge to $N_\mu$.
\end{pf*}

\section{Limit theorems for linear groups}
\label{seccltlin}
%
In this section, we prove the central limit theorem for linear groups
(Theorem~\ref{thcltmat}).
Our main task will be to prove that the norm cocycle \eqref{eqnsigxloggvv}
is centerable.

\subsection{Complete convergence for linear groups}
\label{seccomconlin}
%
In this section, we prove the complete convergence
in the Furstenberg law of large numbers.

Let $\mathbb{K}$ be a local field. The reader who is not familiar with
local fields
may assume $\mathbb{K}=\mathbb{R}$. In general, a local field is a nondiscrete
locally compact field. It is a classical fact that such a field
is a finite extension of either:
\begin{longlist}[(iii)]
\item[(i)]  the field $\mathbb{R}$ of real numbers (in this case, one has
$\mathbb{K}=\mathbb{R}$ or $\mathbb{C}$), or
\item[(ii)]  the field $\mathbb{Q}_p$ of $p$-adic numbers, for some prime number
$p$, or
\item[(iii)]  the field $\mathbb{F}_p((t))$ of Laurent series with
coefficients in the
finite field $\mathbb{F}_p$ of cardinality $p$, for some prime number $p$.
\end{longlist}

Let
$V$ be a finite-dimensional $\mathbb{K}$-vector space.
We fix a basis $e_1,\ldots,e_d$ of $V$ and the following norm on $V$.
For $v=\sum v_ie_i\in V$, we set
$\| v\|= (\sum|v_i|^2)^{{1}/2}$ when $\mathbb{K}=\mathbb{R}$ or
$\mathbb{C}$,
and $\| v\|= \max(|v_i|)$ in the other cases.
We denote by $e^*_1,\ldots,e_d^*$ the dual basis of $V^*$ and
we use the same symbol $\|\cdot\|$ for the norms induced
on the dual space $V^*$, on the space $\operatorname{End}(V)$ of endomorphisms
of $V$,
or on the exterior product $\wedge^2V$, etc.
We equip the projective space $\mathbb{P} (V )$ with the
distance $d$
given, by
\[
d\bigl(x,x'\bigr)=\frac{\llVert v\wedge v'\rrVert }{\llVert v\rrVert  \llVert
v'\rrVert }\qquad \mbox{for $x=\mathbb{K}
v$, $x'=\mathbb{K} v'$ in $\mathbb{P}(V)$.}
\]
For $g$ in $\operatorname{GL}(V)$, we write $N(g):=\max(\|g\|,\|g^{-1}\|)$.

Let $\mu$ be a Borel probability measure on $G:=\operatorname{GL}(V)$ with finite
first moment: $\int_G\log N(g)\,\mathrm{d}\mu(g)<\infty$.
We denote by $\Gamma_\mu$ the subsemigroup of $G$ spanned by the support
of $\mu$,
and by $\lambda_1$ the first Lyapunov exponent of $\mu$,
%
\begin{equation}
\label{eqnlalogmung} \lambda_1:=\lim_{n\rightarrow\infty}
\frac{1}n\int_G\log\| g\|\,\mathrm{d}
\mu^{*n}(g).
\end{equation}

Let $b_1,\ldots, b_n,\ldots$ be random elements of $G$
chosen independently with law~$\mu$.
The Furstenberg law of large numbers
describes the behavior of the random variables
$\log\| b_n\cdots b_1\|$. It is a direct consequence
of the Kingman subadditive ergodic theorem (see, e.g., \cite{Stee}).
It states that,
for $\mu^{\otimes\mathbb{N}^*}$-almost any sequence $(b_1,\ldots
,b_n,\ldots
)$ in $G$,
one has
%
\begin{equation}
\label{eqnloggng1} \lim_{n\rightarrow\infty}\frac{1}n\log\|
b_n\cdots b_1\| =\lambda_1.
\end{equation}

The following Proposition~\ref{procomconlin} is an analogue of
the Baum--Katz theorem for the Furstenberg law of large numbers.
For $p=2$, it says that when the second moment of $\mu$ is finite,
this sequence \eqref{eqnloggng1} converges completely.

\begin{Prop}
\label{procomconlin}
Let $p>1$ and $V=\mathbb{K}^d$.
Let $\mu$ be a Borel\break probability measure on the group $G:=\operatorname{GL}(V)$,
such that the $p$th-moment $\int_G (\log N( g))^p\,\mathrm
{d}\mu(g)$ is finite.
Then, for every $\varepsilon>0$, there exist constants
$C_n=C_n(p,\varepsilon,\mu)$
such that
$\sum_{n\geq1} n^{p-2} C_n<\infty$ and
%
\begin{equation}
\label{eqncomconlin0} \mu^{*n}\bigl(\bigl\{g\in G \mbox{ such that } \bigl|\log\| g
\|-n\lambda_1\bigr|\geq\varepsilon n\bigr\}\bigr)\leq C_n.
\end{equation}
Moreover, if $\Gamma_\mu$
acts irreducibly on $V$,
for any $v$ in $V\setminus\{ 0\}$, one has
%
\begin{equation}
\label{eqncomconlin1} \mu^{*n}\biggl(\biggl\{g\in G \mbox{ such that } \biggl|\log
\frac{\| g v\|}{\|v\|}-n\lambda_1\biggr|\geq\varepsilon n\biggr\}\biggr)\leq
C_n.
\end{equation}
\end{Prop}

\begin{pf}
We first prove the claim \eqref{eqncomconlin0}.
We fix $\varepsilon>0$.
We will apply Proposition~\ref{prollncoc} to the group
$G=\operatorname{GL}(V)$ acting on the projective space $X=\mathbb{P}(V)$ and
to the
norm cocycle
\[
\sigma \dvtx G\times X\rightarrow\mathbb{R};\qquad (g,\mathbb{K} v)\mapsto \log
\frac{\| gv\|}{\|v\|}
\]
for which the function $\sigma_{\sup}$ is $\mathrm{L}^p$-integrable.
According to Furstenberg--Kifer and Hennion theorem
in \cite{FuKi}, Theorem~3.9 and 3.10,  and \cite{Hen84}, Theorem~1 and Corollary~2
(see also \cite{BQRW}, Chapter~3),
the Lyapunov exponent $\lambda_1$ is the upper average of~$\sigma$,
that is,
\[
\lambda_1=\sup_{\nu} \int_{G\times X}
\sigma(g,x)\,\mathrm{d}\mu (g)\,\mathrm{d}\nu(x),
\]
and there exists a unique $\Gamma_\mu$-invariant vector subspace
$V'\subset V$
such that, on one hand, the first Lyapunov exponent
$\lambda'_1$ of the image $\mu'$ of $\mu$ in $\operatorname{GL}(V')$
is strictly smaller than $\lambda_1$,
and, on the other hand,
the image $\mu''$ of $\mu$ in $\operatorname{GL}(V'')$ with $V''=V/V'$
has exponent $\lambda_1$ and the cocycle
$\sigma''\dvtx \operatorname{GL}(V'')\times\mathbb{P}(V'')\rightarrow\mathbb
{R};  (g,\mathbb{K} v)\mapsto
\log\frac{\| gv\|}{\|v\|}$
has\vspace*{1pt} unique average $\lambda_1$.

Since $\lambda_1$ is the upper average of $\sigma$, by Proposition~\ref{prollncoc},
there exist constants $C_n=C_n(p,\varepsilon,\mu)$ such that
$\sum_{n\geq1} n^{p-2} C_n<\infty$ and,
for all $v$ in $V\setminus\{ 0\}$ and $n\geq1$,
%
\begin{equation}
\label{eqncomconlin0M} \mu^{*n}\biggl(\biggl\{g\in G\Big\vert \log
\frac{\| g v\|}{\| v\|}-n\lambda_1\geq\varepsilon n\biggr\}\biggr)\leq
C_n.
\end{equation}
Since $\lambda_1$ is the unique average of $\sigma''$,
using again Proposition~\ref{prollncoc}, one can choose $C_n$ such that,
for all $v''$ in $V''\setminus\{ 0\}$ and $n\geq1$,
%
\begin{equation}
\label{eqncomconlin0m} \mu^{*n}\biggl(\biggl\{g\in G\Big\vert \log
\frac{\| g v''\|}{\| v''\|}-n\lambda_1\notin[-\varepsilon n,\varepsilon n]
\biggr\}\biggr)\leq C_n,
\end{equation}
where, as usual, the norm in the quotient space $V''$ is defined by the equality
$\|v''\|=\inf\{\| v\|\vert  v\in v'' + V'\}$.

The claim \eqref{eqncomconlin0},
with a different constant $C_n$, follows from a combination of
the claim \eqref{eqncomconlin0M}
applied to a basis $v_1,\ldots,v_d$ of $V$
and from the claim \eqref{eqncomconlin0m}
applied to a nonzero vector $v''$ in $V''$.
One just has to notice that there exists a positive constant $M$
such that one has
\[
\log\frac{\| g v''\|}{\|v''\|} \leq\log\| g \|\leq \max_{1\leq i\leq d}\log
\frac{\| g v_i\|}{\|v_i\|}+ M,
\]
for all $g$ in $\operatorname{GL}(V)$ preserving $V'$.

The claim \eqref{eqncomconlin1} follows from \eqref{eqncomconlin0m}, since
when the action of $\Gamma_\mu$ on $V$ is irreducible, one has $V''=V$.
\end{pf}

We denote by $\lambda_2$ the second Lyapunov exponent of $\mu$, that is,
%
\begin{equation}
\label{eqnlalogmung2} \lambda_2 :=\lim_{n\rightarrow\infty}
\frac{1}n\int_G\log\frac{\| \wedge^2 g\|}{\| g\|}\,\mathrm{d}
\mu^{*n}(g).
\end{equation}

\begin{Cor}
\label{corcomconlin}
Assume the same assumptions as in Proposition~\ref
{procomconlin}.
For every $\varepsilon>0$, there exist constants $C_n$
such that
$\sum_{n\geq1} n^{p-2} C_n<\infty$ and
%
\begin{equation}
\label{eqncomconlin4} \mu^{*n}\bigl(\bigl\{g\in G \mbox{ such that } \bigl|\log\|
\wedge^2 g\|-n(\lambda_1 + \lambda_2)\bigr|\geq
\varepsilon n\bigr\} \bigr)\leq C_n.
\end{equation}
\end{Cor}

\begin{pf}
Our statement \eqref{eqncomconlin4} is nothing but \eqref{eqncomconlin0}
applied to
$\wedge^2V$.
\end{pf}

\begin{Rmqs}
An endomorphism $g$ of $V$ is said to be proximal if it admits an
eigenvalue $\lambda$ which has multiplicity one and if all other
eigenvalues of $g$ have modulus ${<}|\lambda|$.
The action of $\Gamma_\mu$ on $V$ is said to be \textit{proximal} if
$\Gamma_\mu
$ contains a proximal endomorphism.
The action of $\Gamma_\mu$ on $V$ is said to be \textit{strongly irreducible}
if no proper finite union of vector subspaces of $V$ is $\Gamma_\mu$-invariant.

According to a result of Furstenberg (see, e.g., \cite{BoLa}), when
$\Gamma
_\mu$ is unbounded,
included in $\operatorname{SL}(V)$ and strongly irreducible in $V$,
the first Lyapunov exponent is positive: $\lambda_1>0$.

According to a result of Guivarc'h in \cite{Gui81}, when the action
of $\Gamma_\mu$ is proximal and strongly irreducible, the first Lyapunov
exponent is simple, that is, one has $\lambda_1>\lambda_2$. We will
use this
fact in the next section.
\end{Rmqs}

\subsection{Log-regularity in projective space}
\label{seclogreqpro}
%
In this section, we prove the log-regularity of the Furstenberg measure
for proximal stronly irreducible representations when
the second moment of $\mu$ is finite.

For any $y=\mathbb{K} f$ in $\mathbb{P}(V^*)$, we set $y^\bot\subset
\mathbb{P}(V)$
for the orthogonal
projective hyperplane:
$y^\bot=\mathbb{P}(\operatorname{Ker} f)$.
For $x=\mathbb{K} v$ in $\mathbb{P}(V)$ and $y=\mathbb{K} f$ in
$\mathbb{P}(V^*)$, we set
\[
\delta(x,y)=\frac{\llvert f(v)\rrvert }{\llVert f\rrVert \llVert
v\rrVert }.
\]
This quantity is also equal to the distance $\delta(x,y)=d(x,y^\bot)$
in $\mathbb{P}(V)$ and to the distance $d(y,x^\bot)$ in $\mathbb{P}(V^*)$.

\begin{Rmq} \label{remuniquestationary1}
Let $\mu$ be a Borel probability measure on $\operatorname{GL}(V)$ such that
$\Gamma_\mu$ is proximal and strongly irreducible on $V$.
Then, due to a result of Furstenberg, $\mu$~admits a unique $\mu
$-stationary Borel probability measure $\nu$ on $\mathbb P(V)$. For
$\beta$-almost any $b$ in $B$,
the sequence of Borel probability measures $(b_1\cdots b_n)_*\nu$
converges to a Dirac measure (see \cite{BoLa}, Section III.4, in the real case
and \cite{BQRW}, Chapter~3, in the general case).
\end{Rmq}

\begin{Prop}
\label{prohausdimpv}
Let $p>1$ and $V=\mathbb{K}^d$.
Let $\mu$ be a Borel probability measure on $G=\operatorname{GL}(V)$
whose $p$th-moment is finite.
Assume that $\Gamma_\mu$ is proximal and strongly irreducible on $V$.
Let $\nu$ be the unique $\mu$-stationary Borel probability measure on
$X=\mathbb{P}(V)$.
Then, for all $y$ in $\mathbb{P}(V^*)$,
%
\begin{equation}
\label{eqnlogregpv} \int_X\bigl|\log\delta(x,y)\bigr|^{p-1}\,
\mathrm{d}\nu(x) \mbox{ is finite},
\end{equation}
and is a continuous function of $y$.
\end{Prop}

\begin{Rmqs}
By a theorem of Guivarc'h in \cite{G90}, when $\mu$
is assumed to have an exponential moment,
the stationary measure $\nu$ is much more regular:
there exists $t>0$ such that
%
\begin{equation}
\label{eqnhausdimpv} \sup_{y\in\mathbb{P}(V^*)} \int_X
\delta(x,y)^{-t}\,\mathrm{d}\nu (x) <\infty.
\end{equation}

The following proof of Proposition~\ref{prohausdimpv}
is similar to our proof in \cite{BQRW} of Guivarc'h theorem, which is
inspired by \cite{BFLM}.

Note that the integral \eqref{eqnlogregpv} may be infinite
when the action of $\Gamma_\mu$ is
assumed to be ``irreducible'' instead of ``strongly irreducible''
(see Example $\ref{exanongau}$).
\end{Rmqs}

Let $K$ be the group of isometries of $(V,\| \cdot\|)$ and
$A^+$ be the semigroup
\[
A^+:=\bigl\{ \operatorname{diag}(a_1,\ldots,a_d)\vert
|a_1|\geq\cdots\geq|a_d|\bigr\}.
\]
For every element $g$ in $\operatorname{GL}(V)$, we choose a decomposition
\[
g=k_{g} a_g\ell_{g} \qquad \mbox{with
$k_{g}$, $\ell_{g}$ in $K$ and $a_g$ in
$A^+$.}
\]
We denote by $x^M_g\in\mathbb{P}(V)$ the
\textit{density point}
of $g$
and by $y^m_g\in\mathbb{P}(V^*)$ the \textit{density point} of $^t  g$,
that is,
\[
x^M_g:=\mathbb{K} k_{g}e_1
\quad \mbox{and}\quad y^m_g:= \mathbb{K}^t
\ell_{g}e_1^*.
\]
We denote by $\gamma_{1}(g)$ the
\textit{first gap}
of $g$, that is, $\gamma_{1}(g):=\frac{\|\wedge^2g\|}{\| g\|^2}$.

\begin{Lem}
\label{lemdegxxmg}
For every $g$ in $\operatorname{GL}(V)$, $x=\mathbb{K} v$ in $\mathbb{P}(V)$
and $y=\mathbb{K} f $
in $\mathbb{P}(V^*)$,
one has:
\begin{longlist}[(iii)]
\item[(i)]
$\delta(x,y^m_g)\leq\frac{\| gv\|}{\| g\|\| v\|}\leq\delta
(x,y^m_g)+ \gamma_{1}(g)$,
\item[(ii)]
$\delta(x^M_g,y)\leq\frac{\| ^tg f\|}{\| g\|\| f\|}\leq\delta
(x^M_g,y)+ \gamma_{1}(g)$,
\item[(iii)]  $d(gx,x^M_g)  \delta(x,y^m_g)\leq\gamma_{1}(g)$.
\end{longlist}
\end{Lem}

\begin{pf}
For all these inequalities, we can assume that $g$ belongs to $A^+$,
that is,
$g=\operatorname{diag}(a_1,\ldots,a_d)$ with $|a_1|\geq\cdots\geq|a_d|$.
We write $v=v_1+v_2$ with
$v_1$ in $\mathbb{K} e_1$ and $v_2$ in the Kernel of $e_1^*$.
One has then
\[
\| g\|=|a_1|,\qquad \gamma_{1}(g)=\frac{|a_2|}{|a_1|}
\quad \mbox{and}\quad \delta\bigl(x,y^m_g\bigr)=
\frac{\|v_1\|}{\|v\|},
\]
\begin{longlist}[(iii)]
\item[(i)]  follows from\vspace*{2pt}
$\| g\| \| v_1\|\leq\| g v\|\leq\| g\| \| v_1\| +|a_2| \|v_2\|$,

\item[(ii)] follows\vspace*{2pt} from \textup{(i)} by replacing $V$ with $V^*$ and $g$ with
$^tg$,

\item[(iii)] follows\vspace*{2pt} from
$d(gx,x^M_g)  \delta(x,y^m_g)=\frac{\| gv_2\|}{\|g v\|}  \frac
{\| v_1\|
}{\|v\|}
\leq\frac{|a_2|}{|a_1|}$.\quad\qed
\end{longlist}
\noqed\end{pf}

\begin{Lem}
\label{lemhausdimpv}
Under the same assumptions as Proposition~\ref{prohausdimpv}, there
exist constants $c>0$, and $C_n>0$ with $\sum_{n\geq
1}n^{p-2}C_n<\infty
$, and
such that, for $n\geq1$,
$x$ in $\mathbb{P}(V)$ and $y$ in $\mathbb{P}(V^*)$, one has
%
\begin{eqnarray}
\label{eqnmundgxmg} \mu^{*n}\bigl(\bigl\{ g\in G\vert d
\bigl(gx,x^M_g\bigr)\geq e^{-cn}\bigr\}\bigr) &
\leq & C_n,
\\
\label{eqnmundexmgy} \mu^{*n}\bigl(\bigl\{ g\in G\vert \delta
\bigl(x^M_g,y\bigr)\leq e^{-cn}\bigr\}\bigr) &
\leq & C_n,
\\
\label{eqnmundegxy} \mu^{*n}\bigl(\bigl\{ g\in G\vert \delta(gx,y)\leq
e^{-cn}\bigr\}\bigr) &\leq & C_n.
\end{eqnarray}
\end{Lem}

\begin{pf}
We set $c=\frac{1}2(\lambda_1-\lambda_2)$ where $\lambda_1$ and
$\lambda_2$
are the first two Lyapunov exponents of $\mu$ (see Section~\ref{seccomconlin}).
According to Guivarc'h theorem in \cite{Gui81},
since the action of $\Gamma_\mu$ is proximal and strongly irreducible,
one has $\lambda_1>\lambda_2$.
According to Proposition~\ref{procomconlin} and its Corollary~\ref
{corcomconlin},
there exist constants $C_n$ such that $\sum_{n\geq1}n^{p-2}C_n<\infty
$ and
such that, for $n\geq1$,
$x=\mathbb{K} v$ in $\mathbb{P}(V)$ and $y=\mathbb{K} f$ in $\mathbb
{P}(V^*)$ with $\|v\|=\|\varphi
\|=1$,
there exist subsets $G_{n,x,y}\subset G$
with $\mu^{*n}(G_{n,x,y})\geq1-C_n$, such that, for $g$ in $G_{n,x,y}$,
the four quantities
\begin{eqnarray*}
&& \biggl\llvert \lambda_1 -\frac{\log\| g\|}{n} \biggr\rrvert,\qquad
\biggl\llvert \lambda_1 -\frac{\log\| gv\|}{n}\biggr\rrvert,\\
&& \biggl\llvert \lambda_1 -\frac{\log\| ^{t}  g\varphi\|}{n}\biggr\rrvert ,\qquad
\biggl\llvert \lambda_1 - \lambda_2 - \frac{\log\gamma_{1}(g)}{n}
\biggr\rrvert
\end{eqnarray*}
are bounded by
$\varepsilon  (\lambda_1-\lambda_2)$ with $\varepsilon=\frac{1}8$.
We will choose $n_0$ large enough, and prove the bounds
\eqref{eqnmundgxmg}, \eqref{eqnmundexmgy} and \eqref{eqnmundegxy}
only for $n\geq n_0$.
We have to check that, for $n\geq n_0$ and $g$ in $G_{n,x,y}$, one has
\[
d\bigl(gx,x^M_g\bigr)\leq e^{-cn},\qquad
\delta\bigl(x_g^M,y\bigr)\geq e^{-cn} \quad
\mbox{and}\quad \delta(gx,y)\geq e^{-cn}.
\]

We first notice that, according to Lemma~\ref{lemdegxxmg}(i), one\vspace*{-2pt} has
\[
\label{eqndexymg0} \delta\bigl(x,y^m_g\bigr)\geq
e^{-2\varepsilon (\lambda_1-\lambda
_2)n}-e^{-(1-\varepsilon)(\lambda_1-\lambda_2)n}
\]
hence, since $n_0$ is arbitrarily large,
%
\begin{equation}
\label{eqndexymg} \delta\bigl(x,y^m_g\bigr)\geq
e^{-3\varepsilon(\lambda_1-\lambda_2)n}.
\end{equation}
But then, using Lemma~\ref{lemdegxxmg}(iii) one gets, for $n_0$ large enough,
%
\begin{equation}
\label{eqndgxxmg0} d\bigl(gx,x^M_g\bigr) \leq
e^{-(1-\varepsilon)(\lambda_1-\lambda_2)n}e^{3\varepsilon (\lambda
_1-\lambda_2)n} = e^{-(1-4\varepsilon)(\lambda_1-\lambda_2)n}.
\end{equation}
This proves \eqref{eqnmundgxmg}.

Applying the same argument as above to $^t  g$ acting on $\mathbb{P}(V^*)$,
the inequality~\eqref{eqndexymg}\vspace*{-4pt} becomes
%
\begin{equation}
\label{eqndxmgy} \delta\bigl(x^M_g,y\bigr)\geq
e^{-3\varepsilon (\lambda_1-\lambda_2)n}.
\end{equation}
This proves \eqref{eqnmundexmgy}.

Hence, combining \eqref{eqndxmgy} with \eqref{eqndgxxmg0}, one gets,
for $n_0$ large\vspace*{-2pt} enough,
\begin{eqnarray*}
\label{eqndgxy0} \delta(gx,y)&\geq & \delta\bigl(x^M_g,y
\bigr) - d\bigl(gx,x^M_g\bigr)
\\[-2pt]
&\geq & e^{-3\varepsilon(\lambda_1-\lambda_2)n}-e^{-(1-4\varepsilon
)(\lambda_1-\lambda_2)n}\geq e^{-4\varepsilon (\lambda_1-\lambda_2)n}.
\end{eqnarray*}
This proves\vspace*{-2pt} \eqref{eqnmundegxy}.
\end{pf}

\begin{pf*}{Proof of Proposition~\ref{prohausdimpv}}
We choose $c$, $C_n$ as in Lemma~\ref{lemhausdimpv}.
We first check that, for $n\geq1$ and $y$ in $\mathbb{P}(V^*)$, one\vspace*{-2pt} has
%
\begin{equation}
\label{eqnnudexyec1n} \nu\bigl(\bigl\{x\in X\vert \delta(x,y)\leq e^{-cn}\bigr
\}\bigr)\leq C_n.
\end{equation}
Indeed, since $\nu=\mu^{*n}*\nu$, one computes using \eqref{eqnmundegxy}
\begin{eqnarray*}
\nu\bigl(\bigl\{x \in X\vert \delta(x,y) \leq e^{-cn}\bigr\}\bigr)& =&
\int_X\mu^{*n}\bigl(\bigl\{ g \in G\vert
\delta(gx,y) \leq e^{-cn}\bigr\}\bigr) \,\mathrm{d}\nu (x)
\\[-2pt]
&\leq & \int_XC_n \,\mathrm{d}\nu(x) =
C_n.
\end{eqnarray*}
Then
cutting the integral \eqref{eqnhausdimpv}
along the subsets $A_{n-1,y}\setminus A_{n,y}$\vspace*{-2pt} where
\[
A_{n,y}:=\bigl\{x\in X\vert \delta(x,y)\leq e^{-cn}\bigr\}
\]
one gets\vspace*{-2pt} the upperbound
\begin{eqnarray*}
\int_X\bigl|\log\delta(x,y)\bigr|^{p-1} \,\mathrm{d}\nu(x)
&\leq & \sum_{n\geq1} c^{p-1}n^{p-1}
\bigl(\nu(A_{n-1,y})-\nu(A_{n,y})\bigr)
\\[-3pt]
&\leq & c^{p-1}+c^{p-1}\sum_{n\geq1}
\bigl((n + 1)^{p-1}-n^{p-1}\bigr) C_n
\\[-3pt]
&\leq & c^{p-1}+(p - 1) 2^p c^{p-1}\sum
_{n\geq1}n^{p-2} C_n,
\end{eqnarray*}
which is finite. This proves \eqref{eqnlogregpv}.

It remains to check the continuity of the function on $\mathbb{P}(V^*)$
\[
\psi^*\dvtx y\mapsto\int_X\bigl|\log\delta(x,y)\bigr|^{p-1}
\,\mathrm{d}\nu(x).
\]
The fact that the above constants $C_n$ do not depend on $y$
tells us that this function $\psi^*$ is a uniform limit
of continuous functions $\psi_n^*$ given by
\[
\psi_n^*\dvtx y\mapsto\int_X\min\bigl(\bigl|\log
\delta(x,y)\bigr|,cn\bigr)^{p-1}\,\mathrm {d}\nu(x).
\]
Hence the function $\psi^*$ is continuous.
\end{pf*}

\subsection{Solving the cohomological equation}
\label{secsolcohequ}
\begin{quotation}
In this section, we prove that the norm cocycle is centerable.
\end{quotation}

We recall that the norm cocycle $\sigma$ on $X=\mathbb{P}(V)$ is the cocycle
\[
\sigma \dvtx \operatorname{GL}(V)\times\mathbb{P}(V)\rightarrow\mathbb{R};
\qquad (g,\mathbb{K} v)\mapsto\log\frac
{\| gv\|}{\|v\|}.
\]

\begin{Prop}
\label{prosolcohequ}
Let $\mu$ be a Borel probability measure on $G=\operatorname{GL}(\mathbb{K}^d)$
whose second moment is finite.
Assume that $\Gamma_{\mu}$ is proximal and strongly irreducible on
$V:=\mathbb{K}^d$.
Then the norm cocycle $\sigma$ on $\mathbb{P}(V)$ is centerable,
that is, satisfies \eqref{eqnspecoc}.
\end{Prop}

\begin{pf}
Let
%
\begin{equation}
\label{eqnphxintsigxg} \varphi \dvtx x\mapsto\int_G\sigma(g,x)\,
\mathrm{d}\mu(g)
\end{equation}
be the expected increase of the cocycle $\sigma$.
We want to find a continuous function $\psi$ on $X$
such that
%
\begin{equation}
\label{eqnphpsmups} \varphi=\psi-P_\mu\psi+\lambda_1,
\end{equation}
where $P_\mu\psi(x)=\int_G\psi(gx)\,\mathrm{d}\mu(g)$, for all
$x$ in $X$,
and where $\lambda_1$ is the first exponent of $\mu$ on $V$.

Let $\check\mu$ be the image of $\mu$
by $g\mapsto g^{-1}$.
We will also denote by $\sigma$ the norm cocycle on $\mathbb
{P}(V^*)$, that is,
the cocycle
\[
\sigma \dvtx \operatorname{GL}(V)\times\mathbb{P}\bigl(V^*\bigr)\rightarrow
\mathbb{R}; \qquad (g,\mathbb{K} f)\mapsto \log\frac{\| f\circ g^{-1}\|}{\|f\|}.
\]
Since the representation of $\Gamma_{\check{\mu}}$ in $V^*$
is also proximal and strongly irreducible, there exists a
unique $\check{\mu}$-stationary probability measure $\nu^*$
on the dual projective space $\mathbb{P}(V^*)$.

Since the second moment of $\mu$ is finite,
according to Proposition~\ref{prohausdimpv},
this measure $\nu^*$ is $\log$-regular.
Hence, the following formula defines a continuous function $\psi$ on $X$:
%
\begin{equation}
\label{eqnppsmps2} \psi(x)=\int_G\log\delta(x,y)\,\mathrm{d}
\nu^*(y),
\end{equation}
where
$\delta(x,y)=\frac{|f(v)|}{\| f\|  \| v\|}$,
for $x=\mathbb{R} v$ in $\mathbb{P}(V)$ and $y=\mathbb{R} f$ in
$\mathbb{P}(V^*)$.

We check the equality,
%
\begin{equation}
\label{eqnsidedesi} \sigma(g,x)=\log\delta\bigl(x,g^{-1}y\bigr)-\log
\delta(gx,y)+\sigma\bigl(g^{-1},y\bigr)
\end{equation}
by computing each side,
\[
\label{eqnsidedesi2} \log\frac{\| gv\|}{\| v\|}= \log\frac{| f(gv)|}{\|f\circ g\| \| v\|}- \log
\frac{| f(gv)|}{\|f\| \| gv\|}+ \log\frac{\| f\circ g\|}{\| f\|}.
\]
Integrating equation \eqref{eqnsidedesi} on $G\times\mathbb{P}(V^*)$
for the measure $ \mathrm{d}\mu(g)\,\mathrm{d}\nu^*(y)$
and using the $\check\mu$-stationarity of $\nu^*$,
one gets \eqref{eqnphpsmups} since $\lambda_1$ is also the first exponent
of $\check\mu$ in~$V^*$.
\end{pf}

\subsection{Central limit theorem for linear groups}
\label{secposvar}
%
The tools we have developed so far allow us to prove not only our
central limit
Theorem \ref{thcltmat} but also a multidimensional version of this theorem.

For $i=1,\ldots,m$, let $\mathbb{K}_i$ be a local field and $V_i$
be a finite-dimensional normed $\mathbb{K}_i$-vector space, and
let $\mu$ be a Borel probability measure on the locally compact group
$G:=\operatorname{GL}(V_1)\times\cdots\times\operatorname{GL}(V_m)$.
We assume that $\Gamma_\mu$ acts strongly irreducibly in each $V_i$.
We consider the compact space
$X=\mathbb{P}(V_1)\times\cdots\times\mathbb{P}(V_m)$.

We denote by $\sigma \dvtx  G\times X\rightarrow\mathbb{R}^m$ the \textit{multinorm cocycle},
that is,
the continuous cocycle given, for $g=(g_1,\ldots,g_m)$ in $G$ and
$x=(\mathbb{K}_1 v_1,\ldots, \mathbb{K}_m v_m)$ in $X$, by
\[
\sigma(g,x):=\biggl(\log\frac{\| g_1v_1\|}{\|v_1\|},\ldots,\log\frac
{\|
g_mv_m\|}{\|v_m\|}\biggr).
\]
We introduce also the function $\kappa\dvtx G\rightarrow\mathbb{R}^m$
given, for $g$ in $G$, by
\[
\kappa(g):=\bigl(\log\| g_1\|,\ldots, \log\| g_m\|\bigr)
\]
and the function $\ell\dvtx G\rightarrow\mathbb{R}^m$ given by
\[
\ell(g):= \lim_{n\rightarrow\infty}\frac{1}n \kappa
\bigl(g^n\bigr),
\]
so that, the $i$th coefficient of $\ell(g)$ is the logarithm of
the spectral radius
of $g_i$.
For $g$ in $G$, we set $N(g)=\sum_{i=1}^mN(g_i)$.

\begin{Rmq}\label{remuniquestationary2}
Let
$\mu$ be a Borel probability measure on the group $\operatorname{GL}(V_1)\times
\cdots\times\operatorname{GL}(V_m)$ such that, for any $1\leq i\leq m$,
$\Gamma_\mu
$ is proximal and strongly irreducible in $V_i$.
By Remark~\ref{remuniquestationary1},
$\mu$ admits a unique $\mu$-stationary Borel probability measure $\nu
_i$ on $\mathbb P(V_i)$ and, for $\beta$-almost any $b$ in\break $B$,
$(b_1\cdots b_n)_*\nu_i$ converges toward a Dirac mass $\delta_{\xi
_i(b)}$ as $n\rightarrow\infty$. One easily shows that this implies
that the image
$\nu$ of $\beta$ by the map
\[
B\rightarrow X;\qquad b\mapsto\bigl(\xi_1(b),\ldots,
\xi_m(b)\bigr)
\]
is the unique $\mu$-stationary Borel probability measure on $X$ (see,
e.g., \cite{BQRW}, Chapter~1).
\end{Rmq}

Here is the multidimensional version of Theorem~\ref{thcltmat}.

\begin{Thm}
\label{thcltmatbis}
Let
$\mu$ be a Borel probability measure on the group $G:=\operatorname{GL}(V_1)\times\cdots\times\operatorname{GL}(V_m)$
such that
$\Gamma_\mu$ acts strongly irreducibly on each $V_i$,
and such that $\int_G (\log N( g))^2\,\mathrm{d}\mu(g)<\infty$.
\begin{longlist}[(a)]
\item[(a)] There exist an element $\lambda$ in $\mathbb{R}^m$, and a
Gaussian law $N_\mu
$ on $\mathbb{R}^m$
such that,
for any bounded continuous
function $F$ on $\mathbb{R}^m$,
one has
%
\begin{equation}
\label{eqnlalogmunv} \int_{G}F \biggl(\frac{\sigma(g,x)-n\lambda}{
\sqrt{n}} \biggr)
\,\mathrm{d}\mu^{*n}(g) \mathop{\longrightarrow}_{n\rightarrow\infty} \int
_{ \mathbb{R}^m}F(t)\,\mathrm{d}N_\mu(t),
\end{equation}
uniformly for $x$ in $X$, and
%
\begin{equation}
\label{eqnlalogmunb} \int_{G}F \biggl(\frac{\kappa( g)-n\lambda}{
\sqrt{n}}
\biggr)\,\mathrm{d}\mu^{*n}(g) \mathop{\longrightarrow}_{n\rightarrow\infty} \int
_{ \mathbb{R}^m}F(t)\,\mathrm{d}N_\mu(t).
\end{equation}
\item[(b)]  When the local fields $\mathbb{K}_i$ are equal to $\mathbb{R}$
and when $\mu$
is supported by
$\operatorname{SL}(V_1)\times\cdots\times\operatorname{SL}(V_m)$,
the support of this Gaussian law $N_\mu$ is the vector subspace $E_\mu$
of $\mathbb{R}^m$
spanned by $\ell(G_\mu)$ where $G_\mu$ is the Zariski closure of
$\Gamma
_\mu$.

\item[(c)]  When $m=1$, $\mathbb{K}_1=\mathbb{R}$ and $\Gamma_\mu$ has
unbounded image in
$\operatorname{PGL}(V_1)$,
the Gaussian law $N_\mu$ is nondegenerate.
\end{longlist}
\end{Thm}

\begin{Rmq}
Point (b)  gives a very practical way to determine the support
of the limit Gaussian law $N_\mu$.
We recall that the \textit{Zariski closure} $G_\mu$ of $\Gamma_\mu$ in $G$
is the smallest subset of $G$ containing $\Gamma_\mu$ which is
defined by polynomial equations. We recall also that the Zariski closure
of a sub-semigroup of $G$ is always a group.
\end{Rmq}

\begin{pf*}{Proof of Theorem~\ref{thcltmatbis}}
(a) We first notice that
equations \eqref{eqnlalogmunv} and \eqref{eqnlalogmunb} are
equivalent since,
for all $\varepsilon>0$, there exists $c>0$ such that,
for all nonzero vector $v_i$ in $V_i$, all $n\geq1$,
\[
\mu^{*n}\bigl(\bigl\{ g\in G \vert c \| g_i\|\leq
\|g_i v_i\|/\| v_i\| \leq\| g_i\|
\bigr\}\bigr) \geq1-\varepsilon
\]
(see, e.g., \cite{BQRWPS}, Lemma~3.2).

First, assume that, for $1\leq i\leq m$, $\Gamma_\mu$ is proximal in $V_i$.
In this case, by Proposition~\ref{prosolcohequ},
in each $V_i$, the norm cocycle is centerable.
Hence, our cocycle $\sigma$ is also centerable. Besides, since by Remark~\ref{remuniquestationary2} $\mu$ admits a unique stationary probability
measure on $X$, $\sigma$ has a unique covariance.
Equation \eqref{eqnlalogmunv} then directly follows
from the central limit Theorem \ref{thcltcoc}.

In general, by Lemma~\ref{lemproxreduction} below,
for any $1\leq i\leq m$, there exists a positive integer $r_i$, a
number $C_i\geq1$ and a finite-dimensional $\mathbb{K}_i$-vector space
$W_i$ equipped with a strongly irreducible and proximal representation
of $\Gamma_\mu$ such that, for any $g$ in $\Gamma_\mu$, one has
\[
C_i^{-1}\|g_i\|^{r_i}_{V_i}
\leq\| g_i \|_{W_i}\leq\|g_i
\|^{r_i}_{V_i}.
\]
Thus, (a) follows from the proximal case applied to the
representations $W_1, \ldots,W_m$.

(b)  We assume now that all the local fields $\mathbb{K}_i$ are equal
to $\mathbb{R}$
and that $\det(g_i)=1$ for all $g$ in $\Gamma_\mu$. We want to describe
the support of the limit Gaussian law $N_\mu$.
Again, by Lemma~\ref{lemproxreduction}, we can assume that all $V_i$'s
are proximal.

According to \cite{Ben97}, Section~4.6, the set $\kappa(\Gamma_\mu
)$ remains
at bounded distance
from the vector space spanned by $\ell(\Gamma_\mu)$.
Hence, the support of $N_\mu$ is included in $E_\mu$.

Conversely, since $\sigma$ is centerable, by \eqref{eqnphmusisimu}, the
covariance $2$-tensor of $N_\mu$
is given by the formula, for all $n\geq1$,
%
\begin{equation}
\Phi_\mu=\frac{1}n\int_{G\times X}\bigl(
\sigma(g,x)-\psi(x)+\psi (gx)-n\lambda\bigr)^2\,\mathrm{d}
\mu^{*n}(g)\,\mathrm{d}\nu(x),
\end{equation}
where $\psi$ is the continuous function in equation \eqref{eqnspecoc}
and $\nu$ is the unique $\mu$-stationary probability measure on $X$.
Let $E_{\Phi_\mu}\subset\mathbb{R}^m$ be the linear span of $\Phi
_\mu$.
For all $g$ in the support of $\mu^{*n}$ and
all $x$ in the support of $\nu$, the element
%
\begin{equation}
\label{eqnsigmaephimu} \sigma(g,x)-\psi(x)+\psi(gx)-n\lambda \qquad \mbox{belongs to }
E_{\Phi_\mu}.
\end{equation}

In particular, let $g$ be an element of $\Gamma_\mu$ which acts in each
$V_i$ as a proximal endomorphism and let
\[
x^+=\bigl(x_1^+,\ldots,x_m^+\bigr),
\]
where, for any $i$, $x_i^+$ is the attractive fixed point of $g$ in
$\mathbb P(V_i)$.
Since $x_i^+$ is an eigenline for $g_i$ whose eigenvalue
has modulus equal to the spectral radius of $g_i$, we have
\[
\sigma\bigl(g,x^+\bigr)=\ell(g).
\]
Since $\Gamma_\mu$ is strongly irreducible in each $V_i$, for any
$x=(x_1,\ldots,x_m)$ in $X$, there exists $h$ in $\Gamma_\mu$ with
$g^nhx\displaystyle{\mathop{\longrightarrow}_{n\rightarrow\infty}} x^+$.
In particular, the support of $\nu$ contains $x^+$, so that,
applying \eqref{eqnsigmaephimu} to the point $x^+$, we get
%
\begin{equation}
\label{eqnellgephimu} \ell(g)\in\mathbb Z\lambda+E_{\Phi_\mu}.
\end{equation}
Now, since the actions on $V_i$ are strongly irreducible, proximal and
volume preserving,
the Zariski closure $G_\mu$ is semisimple. Hence,
by \cite{Ben00}, there exists a subset $\Gamma_1$ of $\Gamma_\mu$ such
that, for any $i$, the elements of $\Gamma_1$ act as proximal
endomorphisms in $V_i$ and that
the closed subgroup of $\mathbb R^m$ spanned by the set $\ell(\Gamma_1)$
in $\mathbb R^m$ is equal to the vector space $E_\mu$ spanned by $\ell
(G_\mu)$. Hence, by \eqref{eqnellgephimu}
this space $E_\mu$ has to be included in $E_{\Phi_\mu}$ and we are done.

(c)
The main difference with point (b) is that the Zariski closure $G_\mu$
may not be semisimple.
The same argument as in (b) tells us that
$\ell([G_\mu,G_\mu])$ is included in $E_{\Phi_\mu}$
and, since the image of $\Gamma_\mu$ in $\operatorname{PGL}(V_1)$ is unbounded,
the group $[G_\mu,G_\mu]$ is also unbounded and one must have
$E_{\Phi
_\mu}=\mathbb{R}$.
\end{pf*}

To deduce the general case in Theorem~\ref{thcltmatbis}(a) from the
one where all the $V_i$ are $\Gamma_\mu$-proximal, we used the following
purely algebraic
lemma.

\begin{Lem} \label{lemproxreduction}
Let $\mathbb{K}$ be a local field, $V$ be a finite-dimensional normed
$\mathbb{K}
$-vector space and
$\Gamma$ be a strongly irreducible sub-semigroup of $\operatorname{GL}(V)$. Let
$r\geq1$ be the proximal dimension of $\Gamma$ in $V$, that is, the least
rank of a nonzero element $\pi$ of the closure
\[
\overline{\mathbb{K}\Gamma}:=\Bigl\{\pi\in\operatorname{End}(V)\big\vert \pi=\lim
_{n\rightarrow\infty}\lambda_ng_n \mbox{ with }
\lambda_n\in\mathbb{K} , g_n\in\Gamma\Bigr\}
\]
and let $W\subset\wedge^rV$ be the subspace spanned by the lines
$\wedge
^r\pi(V)$, where $\pi$ is a rank $r$ element of $\overline{\mathbb
{K}\Gamma}$.
Then:
\begin{longlist}[(a)]
\item[(a)]  $W$ admits a largest proper $\Gamma$-invariant subspace $U$.

\item[(b)]  The action of $\Gamma$ in the quotient $W':=W/U$ is proximal and
strongly irreducible.

\item[(c)]  Moreover, there exists $C\geq1$ such that, for any $g$ in $\Gamma$, one has
%
\begin{equation}
\label{eqncomparenorms} C^{-1}\| g\|^r\leq\|
\wedge^rg\|_{W'}\leq\| g\|^r.
\end{equation}
\end{longlist}
\end{Lem}

\begin{Rmq}
In case $\mathbb{K}$ has characteristic $0$, the action of
$\Gamma$ in
$\wedge^r V$ is semisimple and $W'=W$.
\end{Rmq}

\begin{pf*}{Proof of Lemma~\ref{lemproxreduction}}
(a)
We will prove that $W$ contains a largest proper $\Gamma
$-invariant subspace
and that this space is equal to
\[
U:=\bigcap_\pi\operatorname{Ker}_W
\bigl(\Lambda^r\pi\bigr)\qquad \mbox{where $\pi$ runs among all rank $r$
elements of $\overline{\mathbb{K}\Gamma}$.}
\]
This space $U$ is clearly $\Gamma$-invariant.
We have to check
that the only $\Gamma$-invariant subspace $U_1$ of $W$
which is not included in $U$
is $U_1=W$.
Let $\pi$ be a rank $r$ element of $\overline{\mathbb{K}\Gamma}$
such that $U_1$ is not included in $\operatorname{Ker}(\wedge^r\pi)$. The
endomorphism $\wedge^r\pi$ is proximal and one has
\[
\wedge^r\pi( U_1)\subset U_1.
\]
As $\wedge^r\pi$ has rank one,
one has
\[
\operatorname{Im}\bigl(\wedge^r\pi\bigr)\subset U_1.
\]
Let $\pi'$ be any rank $r$ element of $\overline{\mathbb{K}\Gamma}$.
Since $\Gamma$ is irreducible in $V$, there exists $f$ in $\Gamma$
such that $\pi' f\pi\neq0$. As
$\pi' f\pi$ also belongs to $\overline{\mathbb{K}\Gamma}$, we get
$\operatorname{rk}(\pi' f\pi)=r$
and, since $\wedge^r(\pi' f)$ preserves $U_1$, one has
\[
\operatorname{Im}\bigl(\wedge^r\pi'\bigr)=
\operatorname{Im}\bigl(\wedge^r\bigl(\pi' f\pi\bigr)
\bigr)\subset U_1.
\]
Since this holds for any $\pi'$, by definition of $W$, we get $U_1=W$,
which should be proved.

(b)  The above argument proves also that,
for any rank $r$ element $\pi$ of $\overline{\mathbb{K}\Gamma}$,
one has
%
\begin{equation}
\label{eqnimlar} \operatorname{Im}\bigl(\Lambda^r\pi\bigr)=
\Lambda^r\pi(W) \quad\mbox{and}\quad \operatorname{Im}\bigl(
\Lambda^r\pi\bigr)\not\subset U.
\end{equation}
In particular, the action of $\Gamma$ in the quotient space $W':=W/U$
is proximal.

Let us prove now that the action of $\Gamma$ in $W'$ is strongly
irreducible. Let $U_1,\ldots,U_r$ be subspaces of $W$, all of them
containing $U$, such that $\Gamma$ preserves $U_1\cup\cdots\cup U_r$.
Since $W'$ is $\Gamma$-irreducible, the spaces $U_1,\ldots, U_r$ span
$W$. Let $\Delta\subset\Gamma$ be the sub-semigroup
\[
\Delta:=\{ g\in\Gamma\vert gU_i=U_i \mbox{ for all } 1
\leq i\leq r\}.
\]
There exists a finite subset $F\subset\Gamma$ such that
\[
\Gamma=\Delta F= F\Delta.
\]
In particular, since $\Gamma$ is strongly irreducible in $V$, so is
$\Delta$. Besides, $\Delta$ also has proximal dimension $r$
and, since $\overline{\mathbb{K}\Gamma}=\overline{\mathbb{K}\Delta
} F$, $W$ is also
spanned by
the lines $\operatorname{Im}(\Lambda^r\pi)$ for rank $r$ elements $\pi$ of
$\overline{\mathbb{K}\Delta}$. By applying the first part of the
proof to
$\Delta$,
since the $\Delta$-invariant subspaces $U_i$ span $W$,
one of them is equal to $W$. Therefore, $W'$ is strongly irreducible.

(c)  We want to prove the bounds \eqref{eqncomparenorms}. First, for
$g$ in $\operatorname{GL}(V)$, one has $\|\wedge^r g\|\leq\|g\|^r$. As for $g$
in $\Gamma$, we have $(\wedge^rg)W=W$ and $(\wedge^rg)U=U$, we get
\[
\| \wedge^rg\|_{W'}\leq\| g\|^r.
\]

Assume now there exists a sequence $(g_n)$ in $\Gamma$ with
\[
\| g_n\|^{-r}\| \wedge^rg_n
\|_{W'}\rightarrow 0
\]
and let us reach a contradiction. If $\mathbb{K}$ is $\mathbb R$,
set $\lambda_n=\| g_n\|^{-1}$.
In general, pick $\lambda_n$ in $\mathbb{K}$
such that $\sup_n|\log(|\lambda_n|\| g_n\|)|<\infty$.
After extracting a subsequence, we may assume $\lambda_n
g_n\rightarrow
\pi$,
where $\pi$ is a nonzero element of $\overline{\mathbb{K}\Gamma}$.
In particular, $\pi$ has rank $\geq r$ and we have
$\lambda_n^r\wedge^rg_n\rightarrow\wedge^r\pi$. Thus, since
$\|\lambda_n^r \wedge^rg_n\|_{W'}\rightarrow0$,
we get $\| \wedge^r\pi\|_{W'}=0$, that is,\vspace*{-3pt}
\[
\wedge^r\pi(W)\subset U.
\]
We argue now as in (a).
Let $\pi'$ be a rank $r$ element of $\overline{\mathbb{K}\Gamma}$.
Since $\Gamma$ is irreducible in $V$, there exists $f$ in $\Gamma$
such that $\pi' f\pi\neq0$. Since $\pi' f\pi$ has rank at least $r$,
it has rank exactly $r$, and since $\wedge^r(\pi' f)$ preserves $U$,
one has\vspace*{-3pt}
\[
\operatorname{Im}\bigl(\wedge^r\pi'\bigr)=
\operatorname{Im}\bigl(\wedge^r\bigl(\pi' f\pi\bigr)
\bigr)\subset U.
\]
Since this holds for any $\pi'$, by definition of $W$, we get $U=W$.
This contradiction ends our proof.
\end{pf*}

\begin{Ex}
\label{exanongau}
There exists a finitely supported probability measure $\mu$ on $\operatorname{SL}(\mathbb{R}^d)$
such that $\Gamma_\mu$ is unbounded and acts irreducibly on $\mathbb{R}^d$,
and such that, if we denote by $\lambda_1$ its Lyapunov first exponent,
the random variables
$\frac{\log\| g_n\cdots g_1\|-n\lambda_1}{\sqrt{n}}$
converge in law to a variable which is not Gaussian.
\end{Ex}

Note that, according to Theorem~\ref{thcltmat},
the action of $\Gamma_\mu$ on $\mathbb{R}^d$ cannot be
strongly irreducible.
In our example, the limit law is the law of a random variable $\sup
(\alpha
_1(Z),\ldots,\alpha_m(Z))$
where $Z$ is a $D$-dimensional Gaussian vector and $\alpha_i$ are linear
forms on $\mathbb{R}^D$.
One can prove that this is a general phenomenon.

\begin{pf*}{Proof of Example~\ref{exanongau}}
Set $d=2$ and $\sigma:=
{0\!\!\! \quad-1\choose  \!1 \quad  \ \! \!0}$.
We just choose
$g_i=\sigma^{\varepsilon_i}
{\!e^{x_i} \ \quad \hspace*{3pt} \!\!0\choose   \, 0\!\!\! \quad\hspace*{3pt}  \  e^{-x_i}}$
where $\varepsilon_i$, $x_i$ are independent random variables,
$\varepsilon_i$ takes equiprobable values in $\{0,1\}$ and
$x_i$ are symmetric and real-valued with the same law $\nu\neq\delta_0$.
One can write
$g_n\cdots g_1=
\sigma^{\eta_n}
{\!\!e^{S_n} \hspace*{3pt}\quad \!\! 0\choose  \ 0\!\! \quad \hspace*{3pt} e^{-S_n}}$
with $\eta_n=\varepsilon_1+\cdots+\varepsilon_n $ and
\[
S_n=x_1+(-1)^{\varepsilon_1}x_2+\cdots+
(-1)^{\varepsilon_1+\cdots
+\varepsilon_{n-1}}x_n.
\]
By the classical CLT, the sequence $\frac{S_n}{\sqrt{n}}$ converges in
law to a nondegenerate Gaussian law.
Hence, the sequence $\frac{1}{\sqrt{n}}\log\|g_n\cdots g_1\|=\frac
{|S_n|}{\sqrt{n}}$ converges in law to a non-Gaussian law.
\end{pf*}

\subsection{Central limit theorem for semisimple groups}
\label{seccltsem}
%
In this section, we prove the central limit theorem for
random walks on semisimple Lie groups for a
law $\mu$ whose second moment is finite and
such that $\Gamma_\mu$ is Zariski dense.

This central limit Theorem \ref{thcltmatter} will only be
an intrinsic reformulation
of Theorem~\ref{thcltmatbis}.
Its main interest is that it describes more clearly the support
of the limit Gaussian law.

We first recall the standard notation for semisimple real Lie groups.
Let $G$ be a semisimple connected linear real Lie group,
$\g g$ its Lie algebra,
$K$ a maximal compact subgroup of $G$,
$\g k$ its Lie algebra,
$\g a$ a Cartan subspace of $\g g$ orthogonal to $\g k$ for the Killing form,
and $A$ the subgroup of $G$, $A:=e^{\g a}$.
Let $\g a^+$ be a closed Weyl chamber in $\g a$,
$\g a^{++}$ the interior of $\g a^+$, $A^+=e^{\g a^+}$.
Let $N$ be the corresponding maximal nilpotent connected subgroup
\[
N:=\Bigl\{n\in G\big\vert \forall H\in\g a^{++}, \lim_{t\rightarrow\infty}
e^{-tH}ne^{tH}=1\Bigr\}.
\]
Let $P$ be the corresponding minimal parabolic subgroup of $G$,
that is, $P$ is the normalizer of $N$.
Let $X=G/P$ be the flag variety of $G$.

Using the Iwasawa decomposition $G=KAN$
one defines the \textit{Iwasawa cocycle} $\sigma\dvtx G\times X\rightarrow\g a$:
for $g$ in $G$ and $x$ in $X$,
$\sigma(g,x)$ is the unique element of $\g a$ such that
\[
gk\in Ke^{\sigma(g,x)}N\qquad \mbox{for $x=kP$ with $k$ in $K$}.
\]

Using the Cartan decomposition $G=KA^+K$, one defines
the \textit{Cartan projection} $\kappa\dvtx  G\rightarrow\g a^+$: for $g$ in $G$,
$\kappa(g)$ is the unique element of $\g a^+$ such that\vspace*{-3pt}
\[
g\in Ke^{\kappa(g)}K.
\]

We also define the \textit{Jordan projection} $\ell\dvtx G\rightarrow\g a$ by
\[
\ell(g):=\lim_{n\rightarrow\infty}\frac{1}n \kappa
\bigl(g^n\bigr).
\]

\begin{Exa*}
Before stating the main theorem, let us describe briefly these notions for
$G=\operatorname{SL}(d,\mathbb{R})$. We endow $\mathbb{R}^d$ with the standard
Euclidean
inner product. In this case, one has:
\begin{longlist}[--]
\item[--] $G=\{ g\in \operatorname{End}(\mathbb{R}^d)|\det(g)=1\}$,
$\g g=\{ H\in\operatorname{End}(\mathbb{R}^d)\vert \operatorname{tr}(H)=0\}$,
\item[--] $K=\operatorname{SO}(d,\mathbb{R})=\{g\in G\vert {}^tgg=e\}$,
$\g k=\{ H\in\g g\vert {}^tH+H=0\}$,
\item[--] $\g a=\{ H = \operatorname{diag}(H_1,\ldots, H_d) \in  \g g\}$,
$\g a^+=\{ H\in\g a /  H_1\geq\cdots\geq H_d\}$,
\item[--] $A=\{ a = \operatorname{diag}(a_1,\ldots, a_d) \in  G\vert
a_i>0\}$,
$A^+=\{ a \in  A \vert   a_1 \geq \cdots  \geq  a_d\}$,
\item[--] $N$ is the group of upper triangular matrices with $1$'s on the
diagonal,
\item[--] $P$ is the group of all upper triangular matrices in $G$,
\item[--] $X$ is the set of flags $x=(V_i)_{0\leq i\leq d}$ of $\mathbb{R}^d$,
that is, of increasing sequences of vector subspaces $V_i$
with $\dim V_i=i$.
\item[--] The $i$th coordinate $\sigma_i(g,x)$ of the Iwasawa cocycle
$\sigma(g,x)$ is the logarithm of the norm of the transformation
induced by $g$ between the Euclidean lines $V_i/V_{i-1}\mapsto
gV_i/gV_{i-1}$.
\item[--] The coordinates $\kappa_i(g)$ of the Cartan projection
$\kappa(g)$ are the logarithms of the eigenvalues of $({}^tgg)^{{1}/2}$
in decreasing order.
\item[--] The coordinates $\ell_i(g)$ of the Jordan projection
$\ell(g)$ are the logarithms of the moduli of the eigenvalues of $g$
in decreasing order.
\end{longlist}
\end{Exa*}

\begin{Thm}
\label{thcltmatter}
Let
$\mu$ be a probability measure on the
semisimple connected linear real Lie group $G$.
Assume that
$\Gamma_\mu$ is Zariski dense in $G$,
and that the second moment $\int_G \|\kappa( g)\|^2 \,\mathrm{d}\mu
(g)$ is finite.
Then:
\begin{longlist}[(a)]
\item[(a)]  The Iwasawa cocycle is centerable.
\item[(b)] There exist $\lambda$ in $\g a^{++}$
and a nondegenerate Gaussian law $N_\mu$ on $\g a$
such that,
for any bounded continuous
function $F$ on $\g a$,
one has
%
\begin{equation}
\label{eqnsigxmung} \int_{G}F \biggl(\frac{\sigma(g,x)-n\lambda}{
\sqrt{n}}
\biggr)\,\mathrm{d}\mu^{*n}(g) \mathop{\longrightarrow}_{n\rightarrow\infty} \int
_{ \g a}F(t)\,\mathrm{d}N_\mu(t),
\end{equation}
uniformly for $x$ in $X$, and
%
\begin{equation}
\label{eqnkagmung} \int_{G}F \biggl(\frac{\kappa(g)-n\lambda}{
\sqrt{n}}
\biggr)\,\mathrm{d}\mu^{*n}(g) \mathop{\longrightarrow}_{n\rightarrow\infty} \int
_{ \g a}F(t)\,\mathrm{d}N_\mu(t).
\end{equation}
\end{longlist}
\end{Thm}

We recall that this theorem is due to Goldsheid and Guivarc'h
in \cite{GoGu} and to Guivarc'h in \cite{G08b}
when $\mu$ has a finite exponential moment.

We recall also that the assumption ``$\Gamma_\mu$ is Zariski dense in $G$''
means that, ``every polynomial function on $G$
which is identically zero on $\Gamma_\mu$ is
identically zero on~$G$.''

\begin{pf*}{Proof of Theorem~\protect\ref{thcltmatter}}
(a)  We use the same method as in \cite{Ben97}.
There exists a basis $\chi_1,\ldots,\chi_m$ of $\g a^*$
and finitely many irreducible proximal representations $(V_1,\rho_1)$,
\ldots, $(V_m,\rho_m)$ of $G$
endowed with $K$-invariant norms such that,
for all $g$ in $G$, and $x=hP$ in $X$,
\[
\chi_i\bigl(\kappa(g)\bigr)=\log\| \rho_i(g)\| \quad
\mbox{and}\quad \chi_i\bigl(\sigma(g,x)\bigr)=\log\frac{ \| \rho_i(g)v_i\|}{\| v_i\|},
\]
where $\mathbb{R}v_i$ is the $hPh^{-1}$-invariant line in $V_i$.
It follows then from Theorem~\ref{prosolcohequ} that,
for all $i\leq m$, the
cocycle $\chi_i\circ\sigma$ is centerable. Hence, the Iwasawa
cocycle $\sigma$
is also centerable.

(b)  Using the same argument as in (a), the convergences to a normal law
$N_\mu$
in \eqref{eqnsigxmung} and \eqref{eqnkagmung}
follow from Theorem~\ref{thcltmatbis}.
Theorem~\ref{thcltmatbis} tells us also that the support of $N_\mu$
is the vector subspace of $\g a$ spanned by the set $\ell(G)$.
Since it contains $\g a^+=\ell(A^+)$, this vector subspace is equal to
$\g a$.
\end{pf*}

%





\printaddresses
\end{document}